\documentclass{amsart}
\usepackage{graphicx}
\usepackage{array}
\usepackage{enumerate}

\newtheorem{thm}{Theorem}[section]
\newtheorem{lem}[thm]{Lemma}

\theoremstyle{definition}

\theoremstyle{remark}

\numberwithin{equation}{section}

\begin{document}

\title[A method to compute the truncated theta function]{A nearly-optimal method to compute the truncated theta function, its derivatives, and integrals}
\author[G.A. Hiary]{Ghaith Ayesh Hiary}
\thanks{Preparation of this material is partially supported by the National Science Foundation under agreements No.  DMS-0757627 (FRG grant) and DMS-0635607 (while at the Institute for Advanced Study). This material is based on the author's PhD thesis.}
\address{Pure Mathematics, University of Waterloo, 200 University Ave West, Waterloo, Ontario, Canada, N2L 3G1.}
\email{hiaryg@gmail.com}

\begin{abstract}
A  poly-log time method to compute the truncated theta function, its derivatives, and integrals is presented.  The method is elementary, rigorous, explicit, and suited for computer implementation. We repeatedly apply the Poisson summation formula to the truncated theta function while suitably normalizing the linear and quadratic arguments after each repetition. The method relies on the periodicity of the complex exponential, which enables the suitable normalization of the arguments, and on the self-similarity of the Gaussian, which ensures that we still obtain a truncated theta function after each application of the Poisson summation. In other words, our method relies on modular properties of the theta function. Applications to the numerical computation of the Riemann zeta function and to finding the number of solutions of Waring type Diophantine equations are discussed. 
\end{abstract}

\maketitle

\section{Introduction}
Sums of the form 

\begin{equation} \label{eq:gensum}
\sum_{k=K_1}^{K_2} g(k) \exp(f(k)), \qquad f(x) \in \mathbb{C}[x]\,,\,\,\, g(x)\in \mathbb{C}[x]\,,
\end{equation}

\noindent
arise in areas such as number theory, differential equations, lattice-point problems, optics, and  mathematical physics, among others. For example, one encounters these sums in the context of Diophantine equations and fractional parts of polynomials (\cite{Ko}), solutions of heat and wave equations (\cite{Mu}), counting of integer points lying close to a curve (\cite{Hu}), numerical integration and quadrature formulas (\cite{Ko}), and motion of harmonic oscillators (\cite{Ka}). Due to the importance such sums, there exists an abundance of methods to bound them. For instance, Vinogradov's~\cite{Vi} methods supply such bounds, which along with some involved sieving techniques are used in attacking Goldbach-Waring type problems (see \cite{LWY} for example).

Despite the substantial interest in the sums (\ref{eq:gensum}), comparatively little is known about how to compute them for general values of their arguments. Yet in some settings, it is useful to be able to compute these sums efficiently and accurately. We soon describe two such settings, both of which originate in number theory.

The simplest examples of the exponential sums (\ref{eq:gensum}) occur when $f(x)$ is of degree one, where we obtain the geometric series and its derivatives for which ``closed-form'' formulae are available. The first non-trivial example occurs when $f(x)$ is a quadratic polynomial. In this case, the exponential sum (\ref{eq:gensum}) can be written as a linear combination of exponential sums of the form:
\begin{equation} \label{eq:mainsum}
F(K,j;a,b):=\frac{1}{K^j}\sum_{k=0}^K k^j \exp(2\pi i a k+2\pi i b k^2)\,. 
\end{equation}

Suppose the integer $j$ is not too large. Then in this article, using ideas rooted in analysis, we prove the sum $F(K,j;a,b)$ can be numerically computed to within $\pm\, \epsilon$, for any positive $\epsilon < e^{-1}$, in poly-log time in $K/\epsilon$. The linear and quadratic arguments $a$ and $b$ are any numbers in $[0,1)$, and $j$ is any integer that satisfies \mbox{$0\le j\le O(\log (K/\epsilon)^{\kappa_0})$}, where $\kappa_0$ is any fixed constant.  

More precisely, we obtain  the following upper bound on the number of elementary arithmetic operations (additions, multiplications, evaluations of the logarithm of a positive number, or evaluations of the complex exponential) on numbers of $O((j+1)\log(K/\epsilon))$ bits that our theta algorithm uses.

\begin{thm} \label{thm: mainthm}
There are absolute constants $\kappa_1$, $\kappa_2$, $A_1$, $A_2$, and $A_3$, such that for any positive \mbox{$\epsilon<e^{-1}$}, any integer $K>0$, any integer $j \ge 0$, any \mbox{$a,b \in [0,1)$}, and with  \mbox{$\nu:=\nu(K,j,\epsilon)=(j+1)\log (K/\epsilon)$},  the value of the function $F(K,j;a,b)$ can be computed to within $\pm \, A_1\, \nu ^{\kappa_1} \epsilon$ using \mbox{$\le A_2\, \nu^{\kappa_2}$} arithmetic operations on numbers of $\le A_3\, \nu^2$ bits.
\end{thm}

We remark that a bit complexity bound follows routinely from the arithmetic operations bound in Theorem~\ref{thm: mainthm}. This is because all the numbers that occur in our algorithm have $\le A_3 \, \nu(K,j,\epsilon)^2$ bits.  We do not try to obtain numerical values for the constants $\kappa_1$ and $\kappa_2$ in Theorem~\ref{thm: mainthm}. With some optimization, they probably can be  taken around $3$. Also, in a practical version of the algorithm,  the arithmetic can be performed using substantially fewer than $A_3\,\nu^2$ bits,  and we will likely be able to replace $\nu(K,j,\epsilon)$ with $j+\log (K/\epsilon)$. If we take $\epsilon=K^{-d}$  in the statement of the theorem, then $\nu(K,j,\epsilon) = (d+1)(j+1)\log K$. So the running time of the algorithm becomes $\le A_2\, (d+1)^{\kappa_2}(j+1)^{\kappa_2}(\log K)^{\kappa_2}$ operations. For $d$ and $j$ bounded by any fixed power of $\log K$, this running time is poly-log in $K$. 

We now discuss two applications of the algorithm of Theorem~\ref{thm: mainthm}. For brevity, we  will often  refer to it as the ``theta algorithm'' because $F(K,j;a,b)$ is directly related to the truncated theta function.  

The values of $\zeta(1/2+it)$ on finite intervals are of great interest to number theorists. For example, the numerical verification of the Riemann hypothesis is clearly dependent on such data. There exist several methods to \textit{compute} $\zeta(1/2+it)$, which means methods to obtain the numerical value of $\zeta(1/2+it)$ to within $\pm\,t^{-\lambda}$, for any fixed  $\lambda>0$, and any $t>1$. A well-known approach to computing $\zeta(1/2+it)$ relies on a straightforward application of the Riemann-Siegel formula. The Riemann-Siegel formula has a main sum of length $\lfloor \sqrt{t/(2\pi)}\rfloor$ terms. A simplified version of that formula is:

\begin{equation} \label{eq:rsform}
\zeta(1/2+it)= e^{-i\theta(t)}\,\Re \left(  2\,e^{-i\theta(t)}\sum_{n=1}^{n_1} n^{-1/2} \exp(it \log n)\right) +\Phi_{\lambda}(t)+O(t^{-\lambda})\,,
\end{equation}

\noindent
where $n_1:=\lfloor \sqrt{t/(2\pi)}\rfloor$, and $\theta(t)$ and $\Phi_{\lambda}(t)$ are certain well-understood functions that can be evaluated accurately in $t^{o_{\lambda}(1)}$ operations on numbers of $O_{\lambda}(\log t)$ bits; see~\cite{OS}. (The notation $O_{\lambda}(t)$ or $t^{o_{\lambda}(1)}$ indicates asymptotic constants are taken as $t\to \infty$,  and they depend only on $\lambda$, where we wish to compute $\zeta(1/2+it)$ to within $\pm\, t^{-\lambda}$.)

Our theta algorithm directly leads to a practical method to compute $\zeta(1/2+it)$ to within $\pm\,t^{-\lambda}$ using $t^{1/3+o_{\lambda}(1)}$ operations on numbers of $O_{\lambda}(\log t)$ bits, and requiring $O_{\lambda}(\log t)$ bits of storage. The derivation is explained in a general context in~\cite{Hi} (similar manipulations can also be found in~\cite{Sc} and in~\cite{Ti}, page 99).  A preliminary step in the derivation is to apply appropriate subdivisions and Taylor expansions to the main sum in the Riemann-Siegel formula in order to reduce its computation to that of evaluating, to within $\pm\,t^{-\lambda -1}$, a sum of about $t^{1/3+o_{\lambda}(1)}$ terms of the form $F(K,j;a,b)$, where $K=O(t^{1/6})$, and $0\le j=O_{\lambda}(\log t)$. The power savings now follow because using the theta algorithm, each of the functions $F(K,j;a,b)$ can be evaluated to within $\pm\, t^{-\lambda-2}$ in poly-log time in $t$.

As another simple and direct application of the theta algorithm, we show how to find the number of solutions of a Waring type Diophantine equation. Suppose we wish to find the number of integer solutions to the system:

\begin{equation} \label{eq:ex2}
\sum_{r=1}^s (\alpha_r\, k_r+\beta_r\, k_r^2)-\sum_{r=s+1}^{s+t} (\alpha_r \,k_r+\beta_r\, k_r^2) \equiv 0 \, (\bmod{M})\,, 
\end{equation}

\noindent
where $0\le k_1,\ldots,k_{s+t}\le  K$, and $\alpha_1,\beta_1,\ldots,\alpha_{s+t},\beta_{s+t}$ are some fixed integers. A straightforward calculation reveals that the number of solutions is given by

\begin{equation} \label{eq:ans1}
\frac{1}{M} \sum_{l=0}^{M-1} \left(\prod_{r=1}^s F(K,0; \alpha_r \,l /M,\beta_r\,l/M)\right)\left( \prod_{r=s+1}^{s+t}\overline{F(K,0; \alpha_r\, l /M,\beta_r\, l/M)} \right)\,.
\end{equation}

Using the theta algorithm, the expression (\ref{eq:ans1}) can be evaluated, to the nearest integer say, in \mbox{$M^{1+o(1)} K^{o_{s,t}(1)}$} time. This is already significantly better than a brute-force search. One can also employ the fast Fourier transform to compute (\ref{eq:ans1}) with sufficient accuracy in $M^{1+o(1)} K^{o_{s,t}(1)}+K^{3+o_{s,t}(1)}$ time. But this is less efficient, and it requires temporarily storing large amounts of data. In the special case $M=K$, one can calculate (\ref{eq:ans1}) to the nearest integer  in $M^{1+o(1)} K^{o_{s,t}(1)}$ time using well-known formulae for complete Gauss sums.

In searching for methods to compute $F(K,j;a,b)$, one should make use of the rich structure of the theta function. The theta function, together with variants, occurs frequently in number theory. It is directly related to the zeta function by a Mellin transform, and it has a functional equation as well as other modular properties. So one anticipates that a fast method to compute the truncated theta function will take advantage of this. 

With this in mind, let us motivate the algorithm of Theorem~\ref{thm: mainthm} in the case $j=0$. To this end, recall the following application of Poisson summation due to van der Corput (see~\cite{Ti}, page 74, for a slightly different version). We refer to this application as the {\em{van der Corput iteration}}, although it is not conventionally labelled as such.

\begin{thm}[van der Corput iteration]
Let $f(x)$ be a real function with a continuous and strictly increasing derivative in $s\le x\le t$. Let $f'(s)=\alpha$ and $f'(t)=\beta$. Then

\begin{equation} \label{eq:corput1}
\sum_{s\le k\le t} \exp(2\pi i f(k))=\sum_{\alpha-\eta <v<\beta+\eta}\int_s^t \exp(2\pi i (f(x)-vx))\,dx+\mathcal{R}_{s,t,f}\,, 
\end{equation}

\noindent
where $\mathcal{R}_{s,t,f}= O\left(\log(2+\beta-\alpha)\right)$ for any positive constant $\eta$ less than 1.
\end{thm}

The van der Corput iteration turns a sum of length $t-s$ terms into a sum of length about $\beta-\alpha=f'(t)-f'(s)$ terms, plus a remainder term $\mathcal{R}_{s,t,f}$. In order for this transformation to be a potentially useful computational device, we need $\beta-\alpha\le \tau (t-s)$ for some absolute constant $0\le \tau <1$.  This ensures the new sum is shorter than the original sum. Moreover, we must be able to compute the remainder term $\mathcal{R}_{s,t,f}$, and each of the integrals in the sum over $v$ in (\ref{eq:corput1}), using relatively few operations. For $\eta$ sufficiently small, the latter are precisely the integrals in the Poisson summation formula that {\em{contain a saddle point}}, where an integral is said to contain a saddle point if the function $\frac{d}{dx}\left(f(x)-vx\right)=f'(x)-v$ vanishes for some $x$ in the interval of integration $[s,t]$. So the integrals containing saddle points are determined by:

\begin{equation}\label{eq:spcc}
f'(x)=v\,,\qquad \textrm{for some $x\in [s,t]\,$,}\qquad \Longleftrightarrow \qquad \alpha\le v\le \beta\,.
\end{equation} 

Still, if we simply ensure $\beta-\alpha\le \tau (t-s)$ for some fixed constant $0\le \tau <1$, then the length of the sum over $v$ in (\ref{eq:corput1}) might be of the same order of magnitude as the length of the original sum. For example, if $\tau=1/2$, then we are only guaranteed a cut in the length by 1/2. So the complexity of the problem appears unchanged (in the sense of power-savings). But if we also require the function $\exp(2\pi i f(x))$ to possess some favorable Fourier transform properties that allow us to turn the $v$-terms into ones suited for yet another application of the transformation (\ref{eq:corput1}), then under such hypotheses, one may hope repeated applications of the van der Corput iteration are possible. If they are, one can compute the original sum over $k$ using $\le \log_2 K$ applications of (\ref{eq:corput1}). ($\log_2 x$ is the logarithm of $x$ to base 2.)

These restrictions on $f(x)$ and its Fourier transform are quite stringent. They severely limit the candidate functions for the proposed strategy. Fortunately, the choice \mbox{$f(x)=ax+bx^2$}, which occurs in $F(K,j;a,b)$, is particularly amenable to repeated applications of the van der Corput iteration. Indeed, if we take $s=0$ and $t=K$ in relation (\ref{eq:corput1}), and assume $\lceil a \rceil < \lfloor a+2bK \rfloor$ say, which is frequently the case, then with $f(x)=ax+bx^2$, and for $\eta$ sufficiently small, the transformation (\ref{eq:corput1}) becomes  

\begin{equation} \label{eq:corput2}
\sum_{k=0}^K \exp(2\pi i ak+2\pi i b k^2)=\sum_{v=\lceil a \rceil}^{\lfloor a+2bK \rfloor} \int_0^K \exp(2\pi i a x+2\pi i b x^2-2\pi i vx)\,dx+\,R_1\,, 
\end{equation}

\noindent
where $R_1:=R_1(a,b,K)$. We remark that if the condition $\lceil a \rceil < \lfloor a+2bK \rfloor$ fails, so $\lfloor a+2bK \rfloor \le \lceil a\rceil$, then $b< 1/K$. This means $b$ will be relatively small. For such small $b$, we will use the Euler-Maclaurin summation formula instead of the van der Corput iteration to calculate the sum on the left side in (\ref{eq:corput2});  see \textsection{3.2} for details. That aside, let us write the relation (\ref{eq:corput2}) as 

\begin{equation}
F(K;a,b)=\tilde{F}(\lfloor a+2bK\rfloor;a,b)+R_1\,,
\end{equation}

\noindent
where 

\begin{equation}
\begin{split}
& F(K;a,b) := \sum_{k=0}^K \exp(2\pi i ak+2\pi i b k^2)\,,\\
& \tilde{F}(\lfloor a+2bK\rfloor;a,b) := \sum_{v=\lceil a \rceil}^{\lfloor a+2bK \rfloor} \int_0^K \exp(2\pi i a x+2\pi i b x^2-2\pi i vx)\,dx\,.
\end{split}
\end{equation}

\noindent
We refer to sums of the form $F(K;a,b)$ as quadratic sums.  We recall the following ``self-similarity'' property of the Gaussian:

\begin{equation} \label{eq:ssg}
\int_{-\infty}^{\infty} \exp(\eta t-t^2)\, dt= \sqrt{\pi}\, \exp(\eta ^2/4),\qquad \eta \in \mathbb{C}\,.
\end{equation}

With this setup, we describe the typical iteration of our algorithm. Using the identities in lemma~\ref{lem:l1} in \textsection{4}, as well as conjugation if necessary, it is easily shown the arguments $a$ and $b$ in (\ref{eq:corput2}) can always be normalized so that $a\in [0,1)$ and $b\in [0,1/4]$. The normalization is important, otherwise successive applications of the Poisson summation (in the form of the van der Corput iteration) will essentially cancel each other. Since $b\in [0,1/4]$, the new sum $\tilde{F}(\lfloor a+2bK\rfloor; a,b)$ has length $\lfloor a+2bK\rfloor \le K/2$, which is at most half the length of the original sum. We observe  each term in $\tilde{F}(\lfloor a+2bK\rfloor; a,b)$ is an integral of the form $\int_0^K \exp(2\pi i a x+2\pi i b x^2-2\pi i v x)\,dx$ for some $\lceil a\rceil \le v\le \lfloor a+2bK\rfloor$. And by construction, each such integral contains a saddle-point. We extract the saddle point contribution from each of these integrals. To do so, we first shift the contour of integration to the stationary phase (at an angle of $\pi/4$). Then we complete the domain of integration on both sides to infinity. Last, we use the self-similarity of the Gaussian (\ref{eq:ssg}) to calculate the completed integral explicitly. This yields a new quadratic exponential sum $F(\lfloor a+2bK\rfloor; a/2b,-1/4b)$. Slightly more explicitly, one obtains:

\begin{equation}\label{eq:ver2r}
\tilde{F}(\lfloor a+2bK\rfloor;a,b)=\frac{e^{\pi i /4-\pi i a^2/(2b)}}{\sqrt{2b}}\,F\left(\lfloor a+2bK\rfloor;\frac{a}{2b},-\frac{1}{4b}\right)+R_2\,,
\end{equation}

\noindent
where $R_2:=R_2(a,b,K)$ is a remainder term. It is shown that the original remainder term $R_1$ in (\ref{eq:corput1}), and the new remainder term $R_2$ in (\ref{eq:ver2r}), can both be computed to within $\pm\,\epsilon$ in poly-log time in $K/\epsilon$. Therefore, on repeating the typical iteration at most $\log_2 K$ times, we arrive at a quadratic sum of a small enough length to be evaluated directly. 

In the typical iteration, most of the effort is spent on computing the ``error terms'' $R_1$ and $R_2$. So in order for the overall algorithm to work, it is critical to prove that $R_1$ and $R_2$ can in fact be computed to within $\pm\,\epsilon$ in poly-log time in $K/\epsilon$. This is accomplished in detail in \textsection{3} and \textsection{6}. Briefly though, let us give a heuristic description of why that is. 

The remainder terms $R_1$ and $R_2$ are implicitly defined by relations (\ref{eq:corput2}) and (\ref{eq:ver2r}), respectively.  It is not hard to show these definitions, together with the Poisson summation formula, and the self-similarity of the Gaussian,  imply $R_1$ and $R_2$ must equal the following:

\begin{equation}
\begin{split}
R_1(a,b,K)&= c_{a,b,K}+  PV\sum_{\substack{ v>\lfloor a+2bK\rfloor \\ \textrm{or } v<\lceil a \rceil\,\,\,\,}}\int_0^K \exp(2\pi i a x+2\pi i b x^2-2\pi i v x)\,dx\,,\nonumber\\
R_2(a,b,K)&= d_{a,b}+\sum_{v=\lceil a \rceil}^{\lfloor a+2bK\rfloor} \int_{\substack{x<0 \textrm{ or}\\ x>K\,\,\,\,\,}} \exp(2\pi i a x+2\pi i b x^2-2\pi i v x)\,dx\,,
\end{split}
\end{equation}

\noindent
where $c_{a,b,K}$ and $d_{a,b}$ are certain easily computable quantities, and {\em{PV}} in front of the sum in $R_1$ means the terms of the infinite sum are taken in conjugate pairs. One observes none of the integrals in $R_1$  and $R_2$ contains a saddle point. Because, by construction, $R_1$ consists of precisely the integrals in the Poisson summation formula with no saddle point, while $R_2$ consists of ``complements'' of such integrals, hence, by the monotonicity of $\frac{d}{dx}(ax+bx^2-vx)=a+2bx-v$, they do not contain saddle points themselves. 

The absence of saddle points from the geometric sums $R_1$ and $R_2$ is the reason they do not present any computational difficulty. This is because the absence of saddle points, when combined with suitable applications of Cauchy's theorem, allows for their oscillations to be controlled easily, and in an essentially uniform way. This means the same suitably chosen contour shift can be applied to a large subset of the integrals in $R_1$ (or $R_2$) to ensure rapid exponential decay in the modulus of their integrands. The shifted integrals can thus be truncated quickly, and at a uniform point (after distance about $\log (K/\epsilon)$, where we wish to evaluated $F(K;a,b)$ to within $\pm\, \epsilon$ say). Once truncated, the quadratic part of the integrand, which is $\exp(2\pi i b x^2)$, can be expanded away as a polynomial in $x$ of low degree (since $2\pi bx^2$ no longer gets too large; see \textsection{3} and lemmas~\ref{lem:a1} and~\ref{lem:a2} for the details). One then finds that in computing $R_1$ and $R_2$ the bulk of the computational effort is exerted on integrals of the form

\begin{equation} \label{eq:igf}
h(z,w):=\int_0^1 t^z \exp(w t)\,dt\,,\qquad 0\le z\,,\,\,\,\, z\in \mathbb{Z}\,,\,\,\,\,\Re(w)\le 0\,.
\end{equation}

The function $h(z,w)$ is directly related to the incomplete gamma function. For purposes of our algorithm, the non-negative integer $z$ will be of size $O(\log(K/\epsilon)^{\tilde{\kappa}})$, where $\tilde{\kappa}$ is some absolute constant. In particular, the range of $z$ is quite constrained, which enables a fast evaluation of the integrals (\ref{eq:igf}) via relatively simple methods. But the literature is rich with methods to compute the incomplete gamma function, and consequently $h(z,w)$, for general values of its arguments. These methods are surveyed in great detail by Rubinstein~\cite{Ru}, where they arise in the context of his derivation of a smoothed approximate functional equation for a general class of $L$-functions. 

We further remark that the linear argument $a$, and the quadratic argument $b$, play different roles in the algorithm. Varying the linear argument $a$ corresponds to sliding the sum over $v$ in (\ref{eq:corput2}), whereas varying the quadratic argument $b$ corresponds to compressing, or stretching, the sum. The latter feature greatly accounts for the utility of the van der Corput iteration in the context of the theta algorithm. Also, the role played by the self-similarity of the Gaussian is crucial, because it is the reason we still obtain a quadratic sum after each application of the van der Corput iteration, making its repetition natural to do. 

At the beginning of each iteration, the algorithm normalizes the pair $(a,b)$ to be in $[0,1)\times [0,1/4]$. Afterwards, it computes the remainder terms $R_1(a,b,K)$ and $R_2(a,b,K)$ to within $\pm\,\epsilon$ in poly-log time in $K/\epsilon$. We comment the remainder terms $R_1$ and $R_2$ can still be computed with the same accuracy and efficiency even if we only normalize $(a,b)$ to be in $[0,1)\times [0,1)$. However, the resulting quadratic sum, which is of length $\approx 2bK$, could then be longer the original sum, which is of length $K$. So, although normalizing $(a,b)$ to be in $[0,1)\times [0,1/4]$ is not important to computing the remainder terms $R_1$ and $R_2$ accurately and efficiently in a single iteration, it is important for the recursive step in the algorithm.  

Notice it is not enough to normalize the quadratic argument $b$ so it is in $[0,1/2)$ (this is straightforward to do using the periodicity of the complex exponential and conjugation if necessary). Because if $b\in[0,1/2)$, then $2bK$ could be very close to $K$. So the length of the new sum in the van der Corput iteration, which is $\approx 2bK$, might be very close to the length of the original sum, which is $K$. In particular, we will not have a sufficiently good upper bound on the number of iterations required by theta algorithm. For example, if $b$ starts close to 1/2 mod 1, then its image under the map $b\leftarrow -1/(4b)$, which is the map that occurs in the algorithm, remains close to 1/2 mod 1. The extra ingredient needed to ensure that $b$ is bounded away from $1/2$, that in fact $b\in[0,1/4]$, is the following (easily-provable) identity from lemma~\ref{lem:l1}: 

\begin{equation}
F(K,j;a,b)=F(K, j; a \pm 1/2, b\pm 1/2)=F(K,j; a\mp 1/2, b\pm 1/2)\,. 
\end{equation}

This concludes our sketch of the theta algorithm in the special case $j=0$. For a general $j\ge 0$, the theta algorithm consists of at most $\log_2 K$ iterations. Each iteration acts on $F(K,j;a,b)$ in the following way:

\begin{equation} \label{eq:simit}
F(K,j;a,b) = \sum_{l=0}^j w_{l,j,a,b,K} \,F\left(q_{a,b,K},l;a^*_{a,b},b^*_{a,b}\right) + R_{K,j,a,b} \,,
\end{equation}

\noindent
where $q_{a,b,K}:=\lfloor a+2bK\rfloor$, $a^*_{a,b}:=a/(2b)$, $b^*_{b}:=-1/(4b)$, and the coefficients $w_{l,j,a,b,K}$ are given by formula (\ref{eq:wsj}) in lemma~\ref{lem:a3}.  The remainder term $R_{K,j,a,b}$ is computed to within $\pm\,\epsilon$ in poly-log time in $K/\epsilon$, via the algorithm. A key point is the tuple $(q_{a,b,K},a^*_{a,b},b^*_{b})$ does not depend on $j$. Therefore, the number of new sums $F(.)$ we need to compute in each iteration is always $\le j+1$.  And since the length of each new sum in (\ref{eq:simit}) is $q_{a,b,K}\le (K+1)/2$, the algorithm has to repeat at most $\log_2 K$ times. 

More generally, our method acts on a sum of the form $\sum_{l=0}^j z_l F(K,l;a,b)$ in the following way:

\begin{equation}\label{eq:wwssjj}
\sum_{l=0}^j z_l\, F(K,l;a,b) = \sum_{l=0}^j \tilde{w}_{l,j,a,b,K} \,F\left(q_{a,b,K},l;a^*_{a,b},b^*_{b}\right) + \sum_{l=0}^j R_{K,l,j,a,b}\,,
\end{equation}

\noindent
where $q_{a,b,K}$, $a^*_{a,b}$, and $b^*_{b}$ are the same as in (\ref{eq:simit}), and

\begin{equation}
\tilde{w}_{l,j,a,b,K}:=\sum_{s=l}^j z_{s} w_{l,s,a,b,K}\,.
\end{equation}

\noindent
In \textsection{3}, we show that the coefficients $\tilde{w}_{l,j,a,b,K}$ do not grow too rapidly with each iteration. Specifically, we show that the maximum modulus of $\tilde{w}_{l,j,a,b,K}$ over all iterations of the algorithm is $O(8^j K^2)$, provided the initial coefficients $z_l$ satisfy $\max_{0\le l\le j}|z_l|= O(1)$ say, which is often the case. This bound is rather generous, but it is sharp enough for purposes of our error analysis, and for bounding the number of bits needed by the algorithm to perform its arithmetic operations.  

The presentation is organized as follows. In \textsection{3}, we describe the typical van der Corput iteration. In \textsection{4}, we provide a pseudo-code for the algorithm. In \textsection{5}, it is shown how to compute the related sums 

\begin{equation}\label{eq:gsssm}
G(K,j;a,b):=\sum_{k=1}^K \frac{1}{k^j}\exp (2\pi i a k+2\pi i b k^2)\,,
\end{equation}

\noindent
with a similar complexity and accuracy to $F(K,j;a,b)$. This is done mainly for use in the separate paper~\cite{Hi}.  Finally, in \textsection{6}, we give proofs of various lemmas employed in the previous sections. Section \textsection{6} includes lemmas~\ref{lem:next1} and~\ref{lem:next2}, which are also intended for use in the separate paper~\cite{Hi}. These two lemmas give a complete account of how the theta algorithm behaves, in the case $j=0$, under small perturbations in the linear argument $a$.

\section{Notation}

We let $\lfloor x \rfloor$  denote the largest integer less than or equal to $x$ , $\lceil x \rceil$ denote smallest integer greater than or equal to $x$, $\{x\}$ denote $x-\lfloor x \rfloor$, $\log x$ denote $\log_e x$, and $\exp(x)$ as well as $e^x$ stand for the exponential function (they are used interchangeably). We define $0^0:=1$ whenever it occurs (e.g. in a binomial expansion). For easy reference, we list contours frequently used in later sections:

\begin{displaymath}
\begin{array}{lll}
&C_0:=\{t\,\,|\,\,0\le t < K\}\,, & C_1:=\{K+it\,\,|\,\,0 \le t < K \}\,, \nonumber \\
&C_2:=\{e^{\pi i /4}\,t\,|\, 0 \le t < \sqrt{2}K\}\,, & C_3=\{-it\,|\,0\le t < \infty\}\,, \nonumber \\
&C_4=\{K-it\,|\,0\le t < \infty\}\,, & C_5:= \{e^{\pi i/4}\,t\,|\,-\infty <t < 0\}\,, \nonumber\\
&C_6:= \{e^{\pi i/4}\,t\,|\, \sqrt{2}K \le t < \infty\}\,, & C_7 :=\{e^{-\pi i /4}\,t\,|\, 0 \le t < \sqrt{2}K\}\,, \nonumber\\ 
&C_8:= C_2\cup C_5\cup C_6\,,& C_9 :=\{t\,|\,0 \le t < \infty \}\,.\nonumber
\end{array}
\end{displaymath}

\noindent
Next, define the functions

\begin{equation}
\begin{split}
I_C(K,j;a,b):=&\frac{1}{K^j} \int_{C} t^j \exp(2\pi i a t+2\pi i b t^2)\,dt\,,   \\
J(K,j;M,a,b):=& \frac{1}{K^j} \int_0^{K} t^j \exp \left (-2\pi a  t-2\pi i b  t^2 \right)\frac{1-\exp\left(-2\pi M  t\right)}{\exp(2\pi t)-1}\,dt\,. 
\end{split}
\end{equation}

\noindent
It is convenient to define $\tilde{I}_C(K,j;a,b):=I_C(K,j;ia,-b)$ because it will occur often. Notice $\tilde{I}_C(K,j;a,b)=e^{-\pi i /2} I_{e^{\pi i /2}C}(K,j;a,b)$, so it is essentially a rotation by $\pi/2$. We also define

\begin{displaymath}
\begin{array}{ll}
p\,\,:= p(a)=\lceil a \rceil\,,& \omega_1 :=\omega_1(a)= \lceil a  \rceil -a\,,\\
q\,\,:=q(a,b,K)=\lfloor a+2bK \rfloor\,, & \omega\,\, :=\omega(a,b,K)=\{a+2bK\}\,, \\
p_1:=p_1(a,b,K)=q(a,b,K)-p(a)\,,&  \nu(K,l,\epsilon):=(l+1)\log(K/\epsilon)\,.
\end{array}
\end{displaymath}

\vspace{10pt}

For any $j\ge 0$ and $\epsilon \in (0,e^{-1})$, we say $K$ is $large$ $enough$ if it is satisfies the lower bound $K>\Lambda(K,j,\epsilon)$, where $\Lambda(K,j,\epsilon):=1000\nu(K,j,\epsilon)^6$, and $\nu(K,j,\epsilon):=(j+1)\log (K/\epsilon)$. For example, if $K$ is large enough, then among other consequences, $e^{-K}<(\epsilon/K)^{1000 (j+1)}$. Finally, in the remainder of the paper, any implicit asymptotic constants are absolute, and are applicable as soon as $\epsilon < e^{-1}$, $0\le j$, and $\Lambda(K,j,\epsilon)<K$, unless otherwise is indicated. 

In \textsection{3} and \textsection{4}, we prove Theorem~\ref{thm: mainthm}, which is our main theorem. 

\section{The basic iteration of the algorithm}

Let $j$ be any non-negative integer, $\epsilon$ any number in $(0,e^{-1})$, $K$ any large enough integer, and $(a,b)$ any pair in $[0,1)\times[0,1)$ (the assumption $b\in [0,1/4]$ is not needed in this section, but it is needed in \textsection{4}). Then with $p:=p(a)=\lceil a\rceil$, and $q:=q(a,b,K)=\lfloor a+2bK\rfloor$, either $p<q$ or $q \le p$. The first possibility is the main case, and it is where the algorithm typically spends most of its time. The second possibility is a boundary point that will be handled separately using the Euler-Maclaurin summation.

\subsection{Main case: $p<q$}
Let $p:=p(K,a,b)=\lceil a\rceil -a$, and $q:=q(K,a,b)=\lfloor a+2bK\rfloor$.  Assume $p<q$. By the Poisson summation formula:

\begin{equation}\label{eq:firstsum1}
F(K,j;a,b)=c_{bd}+ PV\sum_{m=-\infty}^{\infty} I_{C_0}(K,j;a-m,b)\,, 
\end{equation}
\noindent
where $\delta_j$ is Kronecker's delta, and $c_{bd}:=c_{bd}(a,b,j,K)=\frac{1}{2}\,\delta_j+\frac{1}{2}\,e^{2\pi i a K+2\pi i b K^2}$ is a boundary term. The notation $PV$ in (\ref{eq:firstsum1}) stands for principal value, so terms are taken in conjugate pairs. Define, 

\begin{equation}
\begin{split}
S_1(K,j;a,b) &:=\sum_{m=p}^{q}I_{C_0}(K,j;a-m,b)\,,  \\
S_2(K,j;a,b) &:= PV \sum_{m \notin [p,q]} I_{C_0}(K,j;a-m,b)\,.
\end{split}
\end{equation} 

\noindent
Therefore,

\begin{equation}\label{eq:fsa}
F(K,j;a,b)=c_{bd}+S_1(K,j;a,b)+S_2(K,j;a,b)\,. 
\end{equation}

Since the boundary term $c_{bd}$ in (\ref{eq:fsa}) can be computed in a constant number of operations on numbers of $O(\log K)$ bits, then it is enough to show how to deal with $S_1(K,j;a,b)$ and $S_2(K,j;a,b)$. We remark the sum $S_1(K,j;a,b)$ corresponds to the terms in the Poisson summation formula that contain a saddle point, and $S_2(K,j;a,b)$ corresponds to the terms that do not contain a saddle point. The plan is to extract the saddle point contributions from $S_1(K,j;a,b)$, which will yield a new (shorter) quadratic exponential sum, plus a remainder term (involving no saddle points) that we will show is computable to within $\pm\,\epsilon$ in poly-log time in $K/\epsilon$. As for $S_2(K,j;a,b)$, whose terms do not contain saddle-points and hence will not contribute to the new quadratic sum, we will show it too can be computed in a similar amount of time and accuracy. 

\subsubsection{The sum $S_1(K,j;a,b)$}
By definition:

\begin{equation}
S_1(K,j;a,b)=\sum_{m=p}^q I_{C_0}(K,j;a-m,b)\,,
\end{equation}

\noindent
where $p=\lceil a\rceil$, $q=\lfloor a+2bK\rfloor$, $C_0:=\{t\,|\, 0\le t<K\}$, and

\begin{equation}\label{eq:ic0yy}
\begin{split}
I_{C_0}(K,j;a-m,b)&:=\frac{1}{K^j}\int_{C_0} t^j \exp(2\pi i (a-m) t+2\pi i b t^2)\,dt\\
&=\frac{1}{K^j}\int_0^K t^j \exp(2\pi i (a-m) t+2\pi i b t^2)\,dt\,.
\end{split} 
\end{equation}

\noindent
The integral $I_{C_0}(K,j;a-m,b)$ contains a saddle-point when $\frac{d}{dt}\,[(a-m)t-bt^2]$ vanishes for some $0\le t\le K$, which is the interval of integration in (\ref{eq:ic0yy}). This occurs precisely when 

\begin{equation}\label{eq:rangems}
0\le (m-a)/(2b) \le K\,\quad \Longleftrightarrow\quad\,  a\le m\le a+2bK\,.
\end{equation}

\noindent
Since (\ref{eq:rangems}) is exactly the range of summation in the definition of $S_1(K,j;a,b)$, then each integral there contains a saddle-point. As stated earlier, we plan to extract the saddle-point contributions from these integrals, which will produce a new shorter quadratic exponential sum of length $\le q+1$ terms. 

To this end, define the contours $C_1:=\{K+it\, | \, 0\le t<K\}$, and $C_2:=\{e^{\pi i /4}\,t\,|\,0\le t<\sqrt{K}\}$.  So $C_1$ and $C_2$ are the two other sides of a right-angle triangle with base $C_0$. By Cauchy's theorem,

\begin{equation} \label{eq:gh1}
S_1(K,j;a,b)=\sum_{m=p}^q I_{C_2}(K,j;a-m,b)-\sum_{m=p}^q I_{C_1}(K,j;a-m,b)\,.
\end{equation}

We first consider the sum $\sum_{m=p}^q I_{C_1}(K,j;a-m,b)$ in (\ref{eq:gh1}). Let us exclude the term corresponding to $m=q$ in that sum for now as it will require a special treatment.  We apply the change of variable $t\leftarrow K+it$ to each integral $I_{C_1}(K,j;a-m,b)$, followed by interchanging the order of summation, then executing the resulting geometric sum,  to obtain

\begin{equation}\label{eq:cfff101}
\begin{split}
\sum_{m=p}^{q-1} I_{C_1}(K,j;a-m,b)=c_1 \sum_{l=0}^j i^l\,\binom{j}{l} \frac{1}{K^l}\int_0^K t^l&\exp\left(-2\pi \omega t-2\pi i b t^2\right) \\
&\times \,\,\frac{1-\exp(-2\pi p_1 t)}{\exp(2\pi t)-1}\,dt\,,
\end{split}
\end{equation}

\noindent
where $\omega = \{a+2bK\}$, $p_1=q-p$, and $c_1:=c_1(a,b,K)=ie^{2\pi i a K + 2\pi i  b K^2}$.  For any integer $M\ge 0$, define:

\begin{equation}
J(K,l;M, a,b):=\frac{1}{K^l}\int_0^K t^l \exp(-2\pi a t-2\pi i b t^2)\,\frac{1-\exp(-2\pi M t)}{\exp(2\pi t)-1}\,dt\,.
\end{equation}

\noindent
Then (\ref{eq:cfff101}) can be expressed as:

\begin{equation}\label{eq:cf1055}
\sum_{m=p}^{q-1} I_{C_1}(K,j;a-m,b)= c_1 \sum_{l=0}^j i^l \binom{j}{l} J(K,l;p_1,\omega,b)\,. 
\end{equation}

The integrand on the right side of (\ref{eq:cfff101}) declines at least like $e^{-2\pi t}$ with $t$ . The rapid decline permits the interval of integration to be truncated quickly, which enables an accurate and efficient evaluation of the (\ref{eq:cfff101}), hence of $J(K,l;p_1,\omega,b)$. (Notice if the term $m=q$ were included in the sum (\ref{eq:cfff101}), the integrand will decline only like $e^{-2\pi \omega t}$,  which might not be fast enough if $\omega$ is very close to zero, and it is the reason that term was excluded earlier.) 

Indeed, according to lemma~\ref{lem:a1}, which is proved via this approach, there exist absolute constants $\kappa_3$, $\kappa_4$, $A_4$, $A_5$, and $A_6$ such that the function $J(K,l;p_1,\omega,b)$ can be evaluated (in terms of short exponential sums) to within $\pm\, A_4\,10^{\kappa_4}\,\nu(K,j,\epsilon)^{\kappa_4} 4^{-j} \epsilon$ using $\le A_5 \,10^{\kappa_3}\,\nu(K,j,\epsilon)^{\kappa_3}$ operations on numbers of $\le A_6\, \nu(K,j,\epsilon)^2$ bits. Notice the reason we built in the factor $4^{-j}$ in the accuracy is because each term in (\ref{eq:cf1055}) is multiplied by $\binom{j}{l}\le 2^j$, and there are $j+1 \le 2^j$ terms. But even if we require $J(.)$ to be computed to within $\pm\,2^{-jd_1}\, K^{-d_2}\,\epsilon$ say for any fixed $d_1$ and $d_2$, then the running time will still be polynomial in $\nu(K,j,\epsilon)$.

As for the term $I_{C_1}(K,j;a-q,b)$, which we excluded earlier, it is treated as follows. Using the change of variable $t\leftarrow K+it$, followed by a binomial expansion, we obtain 

\begin{equation}\label{eq:cfff103}
\begin{split}
I_{C_1}(K,j;a-q,b)&=c_1\sum_{l=0}^j i^l\, \binom{j}{l} \frac{1}{K^l} \int_0^K t^l \exp(-2\pi\omega t-2\pi i b t^2)\,dt \\
&= c_1\sum_{l=0}^j i^l \binom{j}{l} \tilde{I}_{C_0}(K,l;\omega,b)\,.
\end{split} 
\end{equation}

\noindent
where, as before, $c_1:=c_1(a,b,K)=ie^{2\pi i a K + 2\pi i  b K^2}$, and 

\begin{equation} \label{eq:ic77}
\tilde{I}_{C_0}(K,l;\omega,b):=\frac{1}{K^l}\int_0^K t^l \exp(-2\pi \omega t -2\pi i b t^2)\,dt\,.
\end{equation}

\noindent 
The integrand in $\tilde{I}_{C_0}(K,l;\omega,b)$  might not experience rapid exponential decline, because $\omega$ could be very close to zero (recall $\omega:=\{a+2bK\}$, which could get arbitrarily close to zero). One overcomes this difficulty by using Cauchy's theorem: let $C_7:=\{ t e^{-\pi i /4}\,|\,0\le t<\sqrt{2} K\}$ and $\overline{C_1}:=\{K-it\,|\, 0\le t<K\}$, so $C_7$ and $\overline{C_1}$ are the two other sides of a right-angle triangle with base $C_0$; one finds 

\begin{equation}\label{eq:aabv}
\sum_{l=0}^j i^l \binom{j}{l} \tilde{I}_{C_0}(K,l;\omega,b)=\sum_{l=0}^j i^l \binom{j}{l} \tilde{I}_{C_7}(K,l;\omega,b) -\sum_{l=0}^j i^l \binom{j}{l} \tilde{I}_{\overline{C_1}}(K,l;\omega,b)\,.
\end{equation}

The point now is that if $b$ is not too small, the functions $I_{C_7}(.)$ and $I_{\overline{C_1}}(.)$ in (\ref{eq:aabv}) experience more rapid exponential decay, making them much easier to evaluate than $\tilde{I}_{C_0}(K,l;\omega,b)$. Specifically, by lemma~\ref{lem:a1}, each of the functions $\tilde{I}_{C_7}(K,l;\omega,b)$ and $\tilde{I}_{\overline{C_1}}(K,l;\omega,b)$ can be evaluated to within $\pm\, A_7\,10^{\kappa_6}\,\nu(K,j,\epsilon)^{\kappa_6} 8^{-j} \epsilon$ using $\le A_8\,10^{\kappa_5} \,\nu(K,j,\epsilon)^{\kappa_5}$ operations on numbers of $\le A_9\, \nu(K,j,\epsilon)^2$ bits, provided $1\le 2bK$.  Since in this subsection it is assumed $p=\lceil a\rceil < q=\lfloor a+2bK \rfloor$, it follows $a +1\le a+2bK$, and so  $1\le 2bK$.  Put together, we have 

\begin{equation}\label{eq:cf105}
\begin{split}
\sum_{m=p}^{q} I_{C_1}(K,j;a-m,b)= &\, c_1 \sum_{l=0}^j i^l \binom{j}{l} J(K,l;p_1,\omega,b)+c_1\sum_{l=0}^j i^l \binom{j}{l} \tilde{I}_{C_7}(K,l;\omega,b)  \\
& -\,c_1\sum_{l=0}^j i^l \binom{j}{l} \tilde{I}_{\overline{C_1}}(K,l;\omega,b)\,,
\end{split}
\end{equation}

\noindent
where we have shown the right of (\ref{eq:cf105}) side can be computed accurately and efficiently enough for purposes of proving Theorem~\ref{thm: mainthm}. 

Having disposed of the sum $\sum I_{C_1}(.)$ in (\ref{eq:gh1}), we now consider the sum $\sum I_{C_2}(.)$ there. Recall $C_2:=\{e^{\pi i /4}\,t\,|\,0\le t<\sqrt{2}K\}$. We ``complete''  $C_2$ to span the full range $(-\infty,\infty)$. This yields $C_2=C_8-C_5-C_6$, where $C_8:=\{e^{\pi i/4}\,t\,|\,-\infty<t<\infty\}$, $C_5:=\{e^{\pi i /4}\,t\,|\,-\infty<t<0\}$, and $C_6:=\{e^{\pi i /4}\,t\,|\, \sqrt{2}K\le t<\infty\}$. The advantage of rewriting $C_2$ in this way is the following. The integrand in: 

\begin{equation}
I_{C_2}(K, j;a-m,b)= \frac{e^{\pi i (j+1)/4}}{K^j}\int_0^{\sqrt{2}K} t^j\exp(2\pi i e^{\pi i /4} (a-m) t-2\pi b t^2)\,dt\,,
\end{equation}

\noindent 
experiences large oscillations that lead to a tremendous amount of cancellation. Consider, for instance, $|e^{2\pi i e^{\pi i/4} (a-m) t -2\pi b t^2}|$ reaches a maximum of $e^{\pi (m-a)^2/(4b)}$ at the point $0\le t= (m-a)/(2\sqrt{2}b) \le \sqrt{2}K$, while, in comparison, the actual value of the integral is typically much smaller in size. This makes $I_{C_2}(.)$ difficult to evaluate numerically for such $m$. On the other hand, $I_{C_8}(.)$, which still involves a tremendous amount of cancellation, can be evaluated at once via formula (\ref{eq:ssg}), which is the self-similarity property of the Gaussian. Moreover, the extra integrals $I_{C_5}(.)$ and $I_{C_6}(.)$, which were needed to complete $I_{C_2}(.)$, can also be evaluated efficiently because when $m\in\{p,\ldots,q\}$ and $t\notin [0,\sqrt{2}K]$, the integrand $e^{2\pi i e^{\pi i/4} (a-m) t -2\pi b t^2}$ declines rapidly (a consequence of the lack a saddle-point there). Explicitly, since $C_2=C_8-C_5-C_6$,  we have 

\begin{equation} \label{eq:cf201}
\begin{split}
\sum_{m=p}^q I_{C_2}(K,j;a-m,b)=&\, \sum_{m=p}^q I_{C_8}(K,j;a-m,b)-\sum_{m=p}^q I_{C_5}(K,j;a-m,b)\\
&-\sum_{m=p}^q I_{C_6}(K,j;a-m,b)
\end{split}
\end{equation}

\noindent
We consider $\sum_{m=p}^q I_{C_5}(.)$ first.  Let us exclude the term corresponding to $m=p$ from the sum for now because it will require a special treatment. By Cauchy's theorem and a straightforward estimate, we obtain

\begin{equation}\label{eq:ltee1}
\sum_{m=p+1}^q I_{C_5}(K,j;a-m,b) = c_2 J\left(K,j;p_1,\omega_1,b\right) +O(e^{-K})\,, 
\end{equation}

\noindent
where $\omega_1:=\lceil a\rceil -a$, and $c_2:=c_2(j)=(-1)^je^{(j+1)\pi i/2}$.  Like before, the integral $J(.)$ in (\ref{eq:ltee1}) is handled by lemmas~\ref{lem:a1}.  To deal with the special term $m=p$, we relate it to the integral $\tilde{I}_{C_7}(.)$: 

\begin{equation}\label{eq:ltee2}
I_{C_5}(K,j;a-p,b)=c_2\tilde{I}_{C_7}(K,j;\omega_1,b)+O(e^{-K})\,. 
\end{equation}

\noindent
And we already know how to handle $\tilde{I}_{C_7}(.)$ via lemma~\ref{lem:a2}. 

We now consider the sum $\sum_{m=p}^q I_{C_6}(.)$ in (\ref{eq:cf201}). It is not hard to see $I_{C_6}(K,j;a-m,b)$ is bounded by $O(e^{-K}/K)$ for each $m=p,\ldots,q-1$; hence, these terms are negligible due to our assumption $K$ is large enough. And when $m=q$, we have

\begin{equation}
I_{C_6}(K,j;a-q,b) =c_3\, 2^{\frac{j+1}{2}} e^{-2\pi \omega K} \sum_{l=0}^j  \binom{j}{l} \tilde{I}_{C_9}(K,l;\omega-i\omega+2bK+i2bK,-2ib)\,, \nonumber
\end{equation}

\noindent
where $c_3:=c_3(a,j,K)=e^{(j+1)\pi i/4+2\pi i a K}$. The integral $\tilde{I}_{C_9}(.)$ above is also handled by lemma~\ref{lem:a2}.  Finally, the sum $\sum_{m=p}^q I_{C_8}(.)$ in (\ref{eq:cf201}) produces the new quadratic exponential sums since by lemma~\ref{lem:a3}, we obtain

\begin{equation} \label{eq:finsum}
\sum_{m=p}^q I_{C_8}(K,j;a-m,b)=\sum_{s=0}^j w_{s,j,a,b,K} F(q,s;a^*,b^*)-\delta_{1-p}\,  w_{0,j,a,b,K}\,, 
\end{equation}

\noindent
where $a^*\equiv a/(2b) (\textrm{mod }1)$, $b^*\equiv -1/(4b) (\textrm{mod }1)$, and as is apparent from formula (\ref{eq:wsj}) in lemma~\ref{lem:a3}, the coefficients $w_{s,j,a,b,K}$ can be computed to within $\pm,8^{-j} K^{-2} \epsilon$ say for all $s=0,1,\ldots,j$ using $\le A_{10}\,j^2$ operations on numbers of $\le A_{11}\,\nu(K,j,\epsilon)^2$ bits, where $A_{10}$ and $A_{11}$ are absolute constants. 

More generally, if we wish to compute a linear combination of quadratic sums $\sum_{l=0}^j z_l\, F(K,l;a,b)$, rather than a single quadratic sum, then instead of (\ref{eq:finsum}), we obtain

\begin{equation}\label{eq:wwssjj2}
\sum_{l=0}^j z_l \sum_{m=p}^q I_{C_8}(K,l;a-m,b) = \sum_{l=0}^j \tilde{w}_{l,j,a,b,K} \,F\left(q,l;a^*,b^*\right)-\delta_{1-p}\,\tilde{w}_{0,j,a,b,K}\,, 
\end{equation}

\noindent
where 

\begin{equation}
\tilde{w}_{l,j,a,b,K}:=\sum_{s=l}^j z_{s}\, w_{l,s,a,b,K}\,.
\end{equation}

\noindent
And we have the bound 

\begin{equation}\label{eq:anneq}
|\tilde{w}_{l,j,a,b,K}|\le \left(\max_{0\le l\le j} |z_l|\right)\,\sum_{s=l}^j |w_{l,s,a,b,K}|\,. 
\end{equation}

We consider the growth in the coefficients $\tilde{w}_{s,j,a,b,K}$ with each iteration. For our purposes, it suffices to bound the maximum modulus of $\tilde{w}_{s,j,a,b,K}$ over a full run of the algorithm. We examine two scenarios. In the first scenario,  $2\nu(K,j,\epsilon)^3 \le 2bK$ say. Here, as a consequence of lemma~\ref{lem:a4}, we have: ($\log_2 K$ denotes the logarithm to the base 2)

\begin{equation}\label{eq:wwss1}
\begin{split}
\sum_{m=s}^j |w_{s,m,a,b,K}|\le & (2b)^{-1/2}\,e^{1+j/2\nu^3} (1+j/\nu^3) \\
\le & (2b)^{-1/2}\, e^{1+1/\log_2K}  (1+1/\log_2K)\,. 
\end{split}
\end{equation}

\noindent
In the second scenario, $2bK < 2\,\nu(K,j,\epsilon)^3$. We observe this can happen only in the last iteration (because then $q=\lfloor a+2bK\rfloor$ is not large enough, which is a boundary point of the algorithm). There are two possibilities: $q\le p$ or \mbox{$p<q< 2\,\nu(K,j,\epsilon)^3+1$}. In the former case, the algorithm concludes via the Euler-Maclaurin summation method of \textsection{3.2}, and not via the van der Corput iteration. In particular, if $q\le p$, we do not reach the right side of (\ref{eq:wwssjj2}) at all. In the latter case (the case $p<q< 2\,\nu(K,j,\epsilon)^3+1$), lemma~\ref{lem:a4} supplies the bound 

\begin{equation}\label{eq:wwss2}
\sum_{m=s}^j |w_{s,m,a,b,K}|\le  (2b)^{-1/2}\,(j+1) 4^{j+2}\,.  
\end{equation}

\noindent
Therefore, by the bounds (\ref{eq:anneq}), (\ref{eq:wwss1}), and (\ref{eq:wwss2}), and taking into account the algorithm involves $\le \log_2 K$ iterations and $1\le 2bK$, it follows that the maximum modulus of the coefficients $\tilde{w}_{s,j,a,b,K}$ that can occur over a full run of the algorithm is 

\begin{equation} \label{eq:wsjbound1}
\le e^{\log_2 K} (2\log K)^2 (j+1) 4^{j+2} \sqrt{K} = O\left(8^j K^2\right)\,.
\end{equation}

In \textsection{4}, we use the bound (\ref{eq:wsjbound1}) to determine by how much $\epsilon$ needs to be adjusted over a full run of the algorithm so that the final output is still accurate to within $\pm\,A_1\,\nu^{\kappa_1}\epsilon$, as claimed in Theorem~\ref{thm: mainthm}.

\subsubsection{The sum $S_2(K,j;a,b)$}
By definition

\begin{equation}\label{eq:defs22}
S_2(K,j;a,b):=\sum_{m=q+1}^M I_{C_0}(K,j;a-m,b)+\sum_{m=-M}^{p-1} I_{C_0}(K,j;a-m,b)\,.
\end{equation}

\noindent
Let us deal with the subsum $\sum_{m=q+1}^M I_{C_0}(K,j;a-m,b)$ first. If $m> q$, it holds

\begin{equation} \label{eq:est64}
\left|I_{C_0-iT}(K,j;a-m,b)\right| \le (2T)^j e^{-2\pi (1-\omega)T}\int_0^Ke^{-4\pi bT(K-t)}\,dt  \to_{T\to \infty} 0\,,
\end{equation}

\noindent
where the fact $1\le 2bK$ was used to ensure $b$ is bounded from below. So by Cauchy's theorem we can replace $C_0$ with the contours $C_3=\{-it\,|\,0\le t < \infty\}$ and $C_4=\{K-it\,|\,0\le t < \infty\}$, which yields: 

\begin{equation} \label{eq:est63}
\sum_{m=q+1}^M I_{C_0}(K,j;a-m,b)=\sum_{m=q+1}^M I_{C_3}(K,j;a-m,b) -\sum_{m=q+1}^M I_{C_4}(K,j;a-m,b)\,,
\end{equation}

\noindent
(We remark if $j=0$, then (\ref{eq:est64}) holds uniformly in $a\in [0,2]$ and integers $m>q$. Therefore, (\ref{eq:est63}) holds for all $a\in [0,2]$ and integers $m>q=\lfloor a+2bK\rfloor$. This observation is used in the proof of lemma~\ref{lem:next1} later.) Now, by a routine calculation

\begin{equation}\label{eq:jjj1}
\sum_{m=q+1}^M I_{C_3}(K,j;a-m,b)= c_4 J(K,j;M-q,2bK-\omega,b)+O(e^{-K})\,,
\end{equation}

\noindent
where $c_4=:c_4(j)=e^{-(j+1)\pi i /2}$. A similar calculation gives

\begin{equation}\label{eq:jjjjj2}
\sum_{m=q+2}^M I_{C_4}(K,j;a-m,b) =c_5 \sum_{l=0}^j (-i)^{l} \binom{j}{l} J(K,l;M-q-1,1-\omega,b) +O(e^{-K})\,, 
\end{equation}

\noindent
where we isolated the term $I_{C_4}(K,j;a-q-1;b)$ since it will require a special treatment, and where $c_5:=c_5(a,b,K)=-ie^{2\pi i aK+2\pi i b K^2}$ (note $c_5=-c_1$, where $c_1$ as in (\ref{eq:cfff101})). Last, in the case of $I_{C_4}(K,j;a-q-1;b)$, we have

\begin{equation}\label{eq:jjjjj3}
I_{C_4}(K,j;a-q-1,b)=c_5\sum_{l=0}^j (-i)^{l} \binom{j}{l}  \tilde{I}_{C_9}(K,l;1-\omega,b)\,. 
\end{equation}

\noindent
where $C_9 :=\{t\,|\,0 \le t < \infty \}$. As before, the integrals $J(.)$ and $\tilde{I}_{C_9}(.)$, which occur in (\ref{eq:jjj1}), (\ref{eq:jjjjj2}) and (\ref{eq:jjjjj3}), can by computed to within $\pm\,\epsilon$ in polynomial time in $\nu(K,j,\epsilon)$ by lemmas~\ref{lem:a1} and~\ref{lem:a2}.

As for the second subsum $\sum_{m=-M}^{p-1}I_{C_0}(K,j;a-m,b)$ in (\ref{eq:defs22}), the situation is analogous. We simply use the conjugates of the contours $C_3$ and $C_4$,  then repeat the previous calculations with appropriate modifications, which results in the integrals:

\begin{equation}\label{eq:jjj2}
\begin{split}
\sum_{m=-M}^{p-2} I_{\overline{C_3}}(K,j;a-m,b) &= c_6 J(K,j;M+p-1,1-\omega_1,b)i+O(e^{-K})\,, \\
\sum_{m=-M}^{p-1} I_{\overline{C_4}}(K,j;a-m,b) &= c_7 \sum_{l=0}^j \binom{j}{l} i^{l} J(K,l;M+p,2bK-\omega_1,b)+O(e^{-K})\,, 
\end{split}
\end{equation}

\noindent
and

\begin{equation}\label{eq:missjj} 
I_{\overline{C_3}}(K,j;a-p+1,b) = c_6 \tilde{I}_{C_9}(K,j;1-\omega_1,b) \,,
\end{equation}

\noindent
where $c_6:=c_6(j)=e^{(j+1)\pi i /2}$, and  $c_7:=c_7(a,b,K)=i e^{2\pi i aK+2\pi i b K^2}$ (note $c_7=c_1$, where $c_1$ occurs in (\ref{eq:cfff101})). Once again, the functions on the right side in (\ref{eq:jjj2}) and (\ref{eq:missjj}) can by computed to within $\pm\,\epsilon$ in polynomial time in $\nu(K,j,\epsilon)$ by lemmas~\ref{lem:a1} and~\ref{lem:a2}. Finally, the sum $PV\sum_{|m|>M} I_{C_0}(.)$ is bounded as follows:  

\begin{equation}
PV \sum_{|m|>M} I_{C_0}(K,j;a-m,b) = \sum_{m>M} \frac{2}{K^j}\int_0^K t^j \exp(2\pi i a t+2\pi i b t^2)\cos(2\pi m t)\,dt\,. \nonumber
\end{equation}

\noindent
Integrating by parts this is equal to

\begin{equation}
\begin{split}
-\sum_{m>M} \left(\frac{j}{\pi m K^j}  \int_0^K (1-\delta_j) t^{j-1} \exp(2\pi i a t+2\pi i b t^2)\sin(2\pi m t) \,dt \right. \\
 +\left. \frac{2 i}{m K^j}\int_0^K t^j\left(a+2bt\right) \exp(2\pi i a t+2\pi i b t^2) \sin(2\pi m t) \,dt\right)\,. 
\end{split}
\end{equation}

\noindent
By the second mean value theorem, we deduce for $M>2K$ that

\begin{equation}\label{eq:pvvf}
PV \sum_{|m|>M} I_{C_0}(K,j;a-m,b) = O\left(\sum_{m>M} \frac{K}{ m(m-K)}\right)=O\left(\frac{K}{M}\right)\,. 
\end{equation}

\noindent
Finally, take $M=\lceil 8^j K^3 e^{\nu(K,j,\epsilon)}\rceil$ to obtain 

\begin{equation}
PV \sum_{|m|>M} I_{C_0}(K,j; a-m,b) = O(8^{-j} K^{-2}(\epsilon/K)^{j+1})\,,
\end{equation}

\noindent
which suffices in light of our earlier bound (\ref{eq:wsjbound1}) on the maximum modulus of the coefficients $\tilde{w}_{s,j,a,b,K}$ after a full run of the algorithm. We remark that one can let $M$ tend to $\infty$ unless $j=0$, in which case, one can still let $M$ tend to $\infty$ provided the various $J(.)$ integrals are paired appropriately. Of course, this is what one should do in a practical implementation of the algorithm (we do not do this here to simplify the presentation).

In summary, we have shown the following: Let 

\begin{displaymath}
\begin{array}{lll}
c_1:=i\,e^{2\pi i a K+2\pi i b K^2}\,, & c_2:=(-1)^j\,e^{(j+1)\pi i /2}\,, & c_3:=e^{(j+1)\pi i /2+2\pi i a K}\,,\\
c_4 := e^{-(j+1)\pi i /2}\,, & c_5:=-i\, e^{2\pi i a K+2\pi i b K^2}\,, & c_6:= e^{(j+1)\pi i /2}\,.
\end{array}
\end{displaymath}

\vspace{10pt}

\noindent
Let $w_{l,j}:=w_{l,j,a,b,K}$ be defined as in (\ref{eq:wsj}), and let 

\begin{equation}
\tilde{c}_{bd}:=\frac{1}{2}e^{2\pi i a K+2\pi i b K^2}+\frac{1}{2}\delta_j-w_{0,j}\,\delta_{1-p}\,,\nonumber
\end{equation} 

\noindent
where $\delta_j$ is Kronecker's delta. Also define

\begin{displaymath}
\begin{array}{lll}
a^*:= a/(2b)\,, & b^* :=-1/(4b)\,,& q:=\lfloor a+2bK\rfloor\,, \\
\omega:=\{a+2bK\}\,,& \omega_1:=\lceil a\rceil -a\,,& p:=\lceil a \rceil\,,\,\,\,\,\, p_1:=q-p\,.
\end{array}
\end{displaymath}

\vspace{10pt}

\noindent
Then, for $p<q$, $0\le j$, $\epsilon \in (0,e^{-1})$, and $K$ large enough, it holds

\begin{equation}\label{eq:sum1}
F(K,j;a,b)= \sum_{l=0}^j w_{l,j}\,F(q,l;a^*,b^*)+\tilde{S_1}(K,j;a,b)+S_2(K,j;a,b)+\tilde{c}_{bd}\,,
\end{equation}

\noindent
where, for some absolute constant $\tilde{\kappa}_1$, we have

\begin{equation}\label{eq:sum2}
\begin{split}
\tilde{S_1}(K,j;a,b)=& -\,c_1 \sum_{l=0}^j i^l\,\binom{j}{l}\,J(K,l;p_1,\omega,b) - c_2\,J(K,j;p_1,\omega_1,b) \\ 
&\,-\, c_1 \sum_{l=0}^j i^l\,\binom{j}{l}\,\tilde{I}_{C_7}(K,l; \omega,b) + c_1 \sum_{l=0}^j i^l\,\binom{j}{l}\,\tilde{I}_{\overline{C_1}}(K,l;\omega,b)  \\
&\,-\, c_3\, 2^{\frac{j+1}{2}}\,e^{-2\pi \omega K}\,\sum_{l=0}^j \tilde{I}_{C_9}(K,l;\omega-i\omega+2bK+i2bK,-2ib)  \\
&\,-\,c_2\,\tilde{I}_{C_7}(K,j;\omega_1,b) + O(\nu(K,j,\epsilon)^{\tilde{\kappa}_1}\,8^{-j}\,K^{-2}\,\epsilon)\,.
\end{split}
\end{equation}

\begin{equation}\label{eq:sum3}
\begin{split}
S_2(K,j;a,b)=& \,-\,c_5\,\sum_{l=0}^j  (-i)^l \,\binom{j}{l}\,J(K,l;M,1-\omega,b) + c_4\,J(K,j;M,2bK-\omega,b) \\
&\,+\,c_5\,\sum_{l=0}^j i^l\,\binom{j}{l}\,J(K,l;M,2bK-\omega_1,b)+ c_6\,J(K,j;M;1-\omega_1,b)  \\
&\,-\,c_5\,\sum_{l=0}^j(-i)^l\,\binom{j}{l}\,\tilde{I}_{C_9}(K,l;1-\omega,b) +\, c_6\,\tilde{I}_{C_9}(K,j;1-\omega_1,b) \\
&\,+ \,O(\nu(K,j,\epsilon)^{\tilde{\kappa}_1}\,8^{-j}\,K^{-2}\,\epsilon)\,.
\end{split}
\end{equation}

\noindent
And we have shown, with the aid of lemmas~\ref{lem:a1} and~\ref{lem:a2}, that each of the functions on the right side of (\ref{eq:sum2}) and (\ref{eq:sum3}) can be computed to within $O(\nu(K,j,\epsilon)^{\tilde{\kappa}_2}\, 8^{-j}\, K^{-2}\,\epsilon)$ using $O(\nu(K,j,\epsilon)^{\tilde{\kappa}_3})$ operations on numbers of $O(\nu(K,j,\epsilon)^2)$ bits, where the constants $\tilde{\kappa}_2$ and $\tilde{\kappa}_3$ are absolute.

\subsection{Boundary case: $q \le p$}
This occurs when $b$ is very small. We tackle it using the Euler-Maclaurin summation. Without loss of generality, one may assume $K$ is a multiple of 8. So we may write:

\begin{equation}\label{eq:emcsub8}
F(K,j;a,b)=e^{2\pi i a K+2\pi i b K^2}+\frac{1}{K^j}\sum_{m=0}^{7}\sum_{k=mK/8}^{(m+1)K/8-1} k^j \exp(2\pi i a k+2\pi i b k^2)\,. \nonumber
\end{equation}

\noindent
It suffices to  deal with each inner sum in (\ref{eq:emcsub8}) since there are only 8 of them. By a binomial expansion, we have 

\begin{equation}
\frac{1}{K^j}\sum_{k=mK/8}^{(m+1)K/8} k^j \exp(2\pi i a k+2\pi i b k^2)=c_{K,m}\, 8^{-j} \sum_{l=0}^j m^{j-l} \binom{j}{l} F(K_1,l;a_{K,m},b)\,, \nonumber 
\end{equation}

\noindent
where $c_{K,m}:=c_{K,m,a,b}$ is a quickly computable constant of modulus 1, $0\le m< 8$, $K_1:=K/8$, and $a_{K,m}:=a_{K,m,a,b}=a+mbK/4$. 
Using the periodicity of the complex exponential, we can normalize $a_{K,m}$ so it satisfies $-1/2\le a_{K,m} \le 1/2$. 
Since by assumption $q\le p$, then $0\le a+2bK < 2$. So $0\le 2bK<2$, which implies $0\le 2bK_1<1/4$. Therefore, $0\le |a_{K,m}|+2bK_1<3/4$. Put together, we may now assume our task is to compute a quadratic sum $F(K,j;a,b)$ with $|a|+|2bK|< 3/4$. To this end, define:

\begin{equation}
f_{K,j,a,b}(t):=\frac{t^j}{K^j}\exp(2\pi i a t+2\pi i b t^2)\,, \nonumber
\end{equation}

\noindent
By lemma~\ref{lem:a5}, we obtain

\begin{displaymath}
\max_{0\le t \le K} |f^{(N)}_{K,j,a,b}(t)| \le \left(\frac{j+N}{K}+2\pi (|a|+|2bK|)\right)^N\,, 
\end{displaymath}

\noindent
where $f^{(N)}_{K,j,a,b}(t)$ denotes the $N^{th}$ derivative with respect to $t$. Applying the Euler-Maclaurin summation formula to

\begin{equation}
F(K,j;a,b)=\frac{1}{K^j}\sum_{k=0}^K k^j \exp(2\pi i a k+2\pi i b k^2)=:\sum_{k=0}^K f_{K,j,a,b}(k)\,,
\end{equation}

\noindent
yields

\begin{equation}\label{eq:ems}
\begin{split}
F(K,j;a,b)=&\,\int_0^K f_{K,j,a,b}(t)\,dt+\sum_{n=0}^{N} \frac{(-1)^n\,B_{n}}{n!}(f_{K,j,a,b}^{(n-1)}(K)-f_{K,j,a,b}^{(n-1)}(0))\\
&+O\left(\frac{1}{N!}\int_0^K |B_N(\{t\})\,f_{K,j,a,b}^{(N)}(t)|\,dt\right)\,.
\end{split}
\end{equation}

\noindent
where $\{t\}$ denotes the fractional part of $t$, $B_n$ are the Bernoulli numbers, and $B_n(t)$ are the Bernoulli polynomials; so $B_0=1$, $B_1=-1/2$, $B_2=1/6$,$\ldots\,\,$, and $B_0(t)=1$, $B_1(t)=t-1/2$, $B_2(t)=t^2-t+1/6,\ldots\,\,$. 

Taking $N=\lceil 2\log (8^j K^3 /\epsilon)/\log (8/7)+1\rceil$ in (\ref{eq:ems}), it follows from known asymptotics for $B_n$ and $B_n(\{t\})$ (see~\cite{Ru} for instance) that

\begin{equation}\label{eq:emcbfin}
O\left(\frac{2}{(2\pi)^{N}} \int_{0}^K |f^{(N)}(t)|\, dt\right) =O\left(2K (7/8)^{-N}\right)= O(8^{-j} K^{-2} \epsilon)\,. 
\end{equation}

\noindent
Given our earlier bound (\ref{eq:wsjbound1}) on the maximum modulus of the coefficients $w_{s,l,a,b,K}$ after a full run of the algorithm, the bound (\ref{eq:emcbfin}) suffices for purposes of the algorithm. 

Last, the correction terms in (\ref{eq:ems}) can be computed quickly because there are only $\le N+1\le 10\nu(K,j,\epsilon)$ of them, and each can be computed to within $\pm\,\epsilon$ using $O(\nu(K,j,\epsilon^2))$ operations on numbers of $O(\nu(K,j,\epsilon)^2)$ bit via the recursion formula for $f^{(n)}_{K,j,a,b}(t)$ provided in the proof of lemma~\ref{lem:a5}. It only remains to evaluate the integral $\int_0^K f_{K,j,a,b}(t)\,dt$ in (\ref{eq:ems}), which is the main term. But this is equal to $I_{C_0}(K,j;a,b)$, which is handled by lemma~\ref{lem:a2}. 

\section{The algorithm for $F(K,j;a,b)$}

We call a real pair $(a,b)$ $normalized$ if  $(a,b) \in [0,1) \times [0,1/4]$. The normalization is important because sums are converted to integrals via Poisson summation. Therefore, different choices of $a$ or $b$ produce different integrals. We remark it is mainly the normalization of quadratic argument $b$ that truly matters. Normalizing $a$ so that it is in the interval $[0,1)$ is not critical to what follows. For example, it suffices to take $a\in [-m,m]$ for a fixed integer $m>0$. To normalize the arguments $a$ and $b$ properly, we use the following lemma:

\begin{lem} \label{lem:l1}
For any integer $K\ge 0$, any integer $j\ge 0$, and any $a,b\in \mathbb{C}$, the function $F(K,j;a,b)$ satisfies the identities

\begin{equation}
\begin{split}
F(K,j;a,b)&=F(K,j;a+1,b)=F(K,j;a,b+1) \\
&= F(K,j;a\pm 1/2,b \pm 1/2)=F(K,j;a\mp 1/2,b\pm 1/2)\,.  
\end{split}
\end{equation}

\end{lem}

\begin{proof}
This follows from the fact $\exp(2\pi i (z+1))=\exp(2\pi i z)$, and the fact $(k^2\pm k)/2\in \mathbb{Z}$ for any $k\in \mathbb{Z}$.
\end{proof}

As a direct application of the identities in lemma~\ref{lem:l1} we obtain a simple procedure such that starting with any real pair $(a,b)$ it produces a normalized pair $(a_0,b_0) \in [0,1) \times [0,1/4]$ satisfying  

\begin{equation}
F(K,j;a,b)=F(K,j;a_0,b_0)\,,\qquad \textrm{or}\,\qquad F(K,j;a,b)=\overline{F(K,j;a_0,b_0)}\,.
\end{equation}

\noindent
Notice the pair $(a_0, b_0)$ is independent of $K$ and $j$.  The normalization procedure is used in the pseudo-code below to compute $\sum_{l=0}^j z_l\, F(K,l;a,b)$. As before, we let $\nu(K,j,\epsilon):=(j+1)\log(K/\epsilon)$, and $\Lambda(K,j,\epsilon):=1000\nu(K,j,\epsilon)^6$.\\

\begin{itemize}\itemsep10pt
\item  INPUT: Numbers $a,b\in [0,1)$, an integer $K>0$, a positive number $\epsilon \in (0, e^{-1})$, an integer $j\ge 0$, and an array of numbers $z_l$, $l=0,\ldots, j$, with $|z_l|\le 1$ say.
\item  OUTPUT: A complex number $\mathcal{S}$ that equals $\sum_{l=0}^j z_l\, F(K,l;a,b)$ to within $\pm \,A_1\,\nu(K,j,\epsilon)^{\kappa_1} \epsilon$, where $A_1$ and $\kappa_1$ are the absolute constants in Theorem~\ref{thm: mainthm}. 
\item  INITIALIZE: Set $\mathcal{S}=0$, $flag=0$, and $counter=0$. It suffices to perform arithmetic using $A_3\,\nu(K,j,\epsilon)^2$ bit numbers where $A_3$ is the absolute constant in Theorem~\ref{thm: mainthm}.
\end{itemize}

\vspace{10pt}

\begin{enumerate} \itemsep10pt
\item Normalize $(a,b) \leftarrow (a_0,b_0)$ using the identities in lemma~\ref{lem:l1}. This costs a constant number of operations on numbers of $A_3\,\nu(K,j,\epsilon)^2$ bits. If conjugation is needed to normalize $(a,b)$,  set $flag\leftarrow 1$ and $z_l \leftarrow \overline{z_l}$. 
\item Let  $p=\lceil a_0 \rceil$, and $q=\lfloor a_0+2b_0K \rfloor$. These numbers can be calculated using a constant number of operations on numbers of $A_3\,\nu(K,j,\epsilon)^2$ bits. 
\item If $K<\Lambda(K,j,\epsilon)$ (a boundary case), evaluate the sum $\sum_{l=0}^j z_l\,F(K,l;a,b)$ directly. This can be done using $\le \tilde{A}_1\,(j+1) \Lambda(K,j,\epsilon)$ operations on number of $A_3\, \nu(K,j,\epsilon)^2$ bits, where $\tilde{A}_1$ is an absolute constant. Store the result in $R[counter]$. If $flag=1$,  set  $R[counter]\leftarrow \overline{R[counter]}$. Go to (9). 
\item If $q \le p$ (a boundary case), apply the Euler-Maclaurin technique of \textsection{3.2} to evaluate the sum to within $\pm\,\tilde{\epsilon}$ where $\tilde{\epsilon}:=8^{-j} K^{-2} \epsilon$. This costs $\le \tilde{A}_2\,\nu(K,j,\tilde{\epsilon})^{\tilde{\kappa}_4}$ operations on numbers of $A_3\,\nu(K,j,\epsilon)^2$ bits, where $\tilde{A}_2$ and $\tilde{\kappa}_4$ are absolute constants. (Notice $\nu(K,j,\tilde{\epsilon})\le  4(j+1)\nu(K,j,\epsilon)$, and so  $\tilde{A}_3 \nu(K,j,\tilde{\epsilon})^{\tilde{\kappa}_4}\le 4^{\tilde{\kappa}_4} \tilde{A}_3\,\nu(K,j,\epsilon)^{2\tilde{\kappa}_4}$.) Store the result in $R[counter]$. If $flag=1$, set $R[counter]\leftarrow \overline{R[counter]}$. \mbox{Go to (9)}.
\item  Apply the algorithm iteration for the case $p<q$. This step requires the calculation of the quantities $q:=\lfloor a_0+2b_0 K\rfloor$, $a^*:=\frac{a_0}{2b_0}$, and $b^*:=-\frac{1}{4b_0}$, all of which can be calculated using a constant number of operations.  We obtain

\begin{equation}
\sum_{l=0}^j z_l\, F(K,l;a,b) = \sum_{l=0}^j \tilde{w}_{l,j,a,b,K} \,F\left(q,l;a^*,b^*\right) + \sum_{l=0}^j R_{K,l,j,a,b}\,.\nonumber
\end{equation}

\noindent
where $\tilde{w}_{l,j,a,b,K}:=\sum_{s=l}^j z_{s} w_{l,s,a,b,K}$. The remainder $\sum_{l=0}^j R_{K,l,j,a,b}$ is computed by the algorithm to within $\pm\,\tilde{A}_4 \,\nu(K,j,\epsilon)^{\tilde{\kappa}_5}\epsilon$ using $\le \tilde{A}_5\,\nu(K,j,\tilde{\epsilon})^{\tilde{\kappa}_6}$ operations on numbers of $A_3\,\nu(K,j,\epsilon)^2$ bits, where $\tilde{A}_4$, $\tilde{A}_5$, $\tilde{\kappa}_5$, and $\tilde{\kappa}_6$, are absolute constants.

\item  Set $R[counter]=\sum_{l=0}^j z_{l}R_l$, $z_{l} \leftarrow \sum_{s=l}^j z_{s,j} w_{l,s,a_0,b_0,K}$, $a \leftarrow a^*$, $b \leftarrow a^*$, $K\leftarrow q$, and $counter \leftarrow counter+1$.
\item  If $flag=1$, set $z_{l}\leftarrow \overline{z_{l}}$, $R[counter] \leftarrow \overline{R[counter]}$, $a \leftarrow -a$, $b \leftarrow -b$, and $flag \leftarrow 0$.
\item  Go to (1).
\item  Set $\mathcal{S} = \sum_{l=0}^{counter} R[l]$. Return $\mathcal{S}$.
\end{enumerate}

\section{The sums $G(K,j;a,b)$}
We show how evaluate the sums $G(K,j;a,b)$ defined in (\ref{eq:gsssm}) to within $\pm\,\epsilon$. Assume $K$ is large enough (i.e. $K>\Lambda(K,j,\epsilon)$), otherwise we can evaluate the sum directly.  Define 

\begin{equation}
\tilde{G}(N,j;a,b) :=\sum_{k=N}^{2N-1} \frac{1}{k^j}\exp (2\pi i a k+2\pi i b k^2)\,. 
\end{equation}

\noindent
It is not too hard to show $G(K,j;a,b)$ can be written as the sum of $O(\log K)$ subsums of the form $\tilde{G}(N,j;a,b)$, with $N<K$,  plus a remainder sum of length $O(\log K)$ terms. So it is enough to show how to compute $\tilde{G}(.)$ to within $\pm\, \epsilon$. Without loss of generality, we may assume $N$ is a multiple of 16, so we may write: 

\begin{equation}
\tilde{G}(N,j;a,b)=\sum_{m=0}^{15}\sum_{k=N_m}^{N_{m+1}-1} \frac{1}{k^j}\exp (2\pi i a k+2\pi i b k^2)\,, 
\end{equation}

\noindent
where $N_m:=N+m N/16$. The inner sum in the last expression is

\begin{equation} \label{eq:gsu1}
\frac{c_{N,m}}{N_m^j}\sum_{l=0}^{\infty} (-1)^l \binom{j+l-1}{j-1}\sum_{k=0}^{N/16-1}  \frac{ k^l}{N_m^l}\exp (2\pi i a_{N,m} k+2\pi i b k^2)\,, 
\end{equation}

\noindent
where $c_{N,m}:=c_{N,m,a,b}$ satisfies $|c_{N,m}|=1$, and $a_{N,m}:=a+2bN_m$. Since $\binom{j+l-1}{j-1} \,k^l/N_m^{l+j} \le 8^{-l}$, we can truncate the sum over $l$ in (\ref{eq:gsu1}) after $\lceil 10\log(K/\epsilon)\rceil$ terms say, which yields a truncation error of $\pm \epsilon/K$ say. Finally, by Theorem~\ref{thm: mainthm}, each inner sum in (\ref{eq:gsu1}) can be computed to within $\pm\, \epsilon/K$, using $\le 2^{\kappa_1} A_2\,\nu(K,j,\epsilon)^{\kappa_1}$ operations on numbers of $\le A_3\,\nu(K,j,\epsilon)^2$ bits.

\section{Auxiliary results}

\begin{lem} \label{lem:a1}
There are absolute constants $\kappa_3$, $\kappa_4$, $A_4$, $A_5$, and $A_6$, such that for any positive $\epsilon <e^{-1}$, any integer $0\le j$, any integer $10\,\nu(K,j,\epsilon)^2<K$ say,  any integer $0<M < e^{10\,\nu(K,j,\epsilon)^2}$ say, any $0\le w < K$ say, and any $0\le b\le 1$, the integral $J(K,j;M,w,b)$ can be evaluated to within $\pm\, A_4\,\nu(K,j,\epsilon)^{\kappa_3} \epsilon$ using $\le A_5\,\nu(K,j,\epsilon)^{\kappa_4}$ operations on numbers of $\le A_6\, \nu(K,j,\epsilon)^2$ bits.
\end{lem}

\begin{proof}
The integrand in $J(K,j;M,w,b)$ declines exponential fast, so the integral can be truncated quickly. Specifically, let  $L:=L(K,j,\epsilon)=\lceil \nu(K,j,\epsilon) \rceil$, then 

\begin{equation}
J(K,j;M,w,b)=\frac{1}{K^j} \int_0^{ L} t^j \exp \left (-2\pi w t-2\pi i b t^2 \right)\frac{1-\exp\left(-2\pi M t\right)}{\exp(2\pi t)-1}\,dt +O(\epsilon)\,. \nonumber
\end{equation}

\noindent
Therefore, in order to evaluate $J(K,j;M,w,b)$ in a time complexity as stated in the lemma, it suffices to deal with the integrals

\begin{equation}
g(j,M,w,b,n):=\frac{1}{K^j} \int_n^{n+1} t^j \exp \left (-2\pi w t-2\pi i b t^2 \right)\frac{1-\exp\left(-2\pi M t\right)}{\exp(2\pi t)-1}\,dt\,,\nonumber 
\end{equation}

\noindent
where $n\in\{0,\ldots,L-1\}$. By the change of variable $t\leftarrow t-n$, followed by Taylor expansions applied to the quadratic factor $e^{-2\pi i b t^2}$, we obtain after some simple estimates that

\begin{equation}
\begin{split}
&g(j,M,w,b,n)=\frac{\exp(-2\pi w n -2\pi i b n^2)}{K^j}\sum_{s=0}^j \binom{j}{s} n^{j-s} \sum_{r=0}^{L} \frac{(-2\pi ib)^r}{r!} \times   \\
& \int_0^{1} t^{s+2r} \exp \left (-2\pi w t- 4\pi i b nt \right) \frac{1-\exp\left(-2\pi M (t+n)\right)}{\exp(2\pi (t+n))-1}\,dt +O(\epsilon\,\log M)\,. \nonumber
\end{split}
\end{equation}

\noindent
Since the last expression is a linear combination of $(L+1)(j+1)\le 10\,\nu(K,j,\epsilon)^2$ terms of the form

\begin{equation} \label{eq:add1}
\int_0^{1} t^\alpha \exp \left (-2\pi w t-4\pi i b nt \right)\frac{1-\exp\left(-2\pi M (t+n)\right)}{\exp(2\pi(t+n))-1}\,dt\,,
\end{equation}

\noindent
for integers $0\le \alpha\le 2L+j$, then our task is reduced to dealing with the integral (\ref{eq:add1}) over that range of $\alpha$. To evaluate this integral, we first unfold the geometric series in the integrand; that is, we write (\ref{eq:add1}) as

\begin{equation} \label{eq:add2}
\sum_{m=1}^{M}\exp(-2\pi mn)\int_0^{1} t^\alpha \exp \left(-2\pi (m+w +2i b n )t \right)\,dt\,. 
\end{equation}

\noindent
(Notice the integrals occurring in (\ref{eq:add2}) are incomplete Gamma functions, which we alluded to earlier in formula (\ref{eq:igf}). Although the methods given in this lemma to evaluate such integrals suffice for complexity bounds, there are other more practical, though more tedious to describe, methods). Define $m_{\alpha,n}:=m_{\alpha,n,w}= \max\{1,\lceil \alpha - w - 2bn\rceil\}$, in particular $\alpha \le m_{\alpha,n}+w+2bn$. We split (\ref{eq:add2}) into two subsums: $\sum_{m_{\alpha,n}\le m\le M}$ and $\sum_{1 \le m< m_{\alpha,n}}$ (the splitting of the sum is because the general function $h(z,w):=\int_0^1 t^z\exp(wt)\,dt$ behaves essentially differently according to whether $|w|<|z|$ or $|z|<|w|$). Each term in the subsum $\sum_{m_{\alpha,n} \le m\le M}$ can be calculated explicitly as

\begin{equation}\label{eq:anintg2}
\begin{split}
\int_0^{1} t^\alpha \exp \left(-2\pi (m+w +2i b n )t \right)\,dt=&-\sum_{v=1}^{\alpha+1} \frac{\alpha !}{(\alpha+1-v)!}\,\frac{\exp(-2\pi m-2\pi w- 4\pi b i n ))}{(2\pi m+2\pi w+4\pi i b n)^v} \\
&+ \frac{\alpha !}{(2\pi m+2\pi w+4\pi i b n)^{\alpha+1}}\,.
\end{split}
\end{equation}

\noindent
So, on interchanging the order of summation, the subsum $\sum_{m_{\alpha,n}\le m\le M}$  is equal to

\begin{equation} \label{eq:aux4}
\begin{split}
-\sum_{v=1}^{\alpha+1} \frac{\alpha !}{(\alpha+1-v)!}\sum_{m=m_{\alpha,n}}^{M} & \exp(-2\pi mn)\,\frac{\exp(-2\pi m-2\pi w- 4\pi b i n ))}{(2\pi m+2\pi w+4\pi i b n)^v}  \\
&+\, \alpha !\sum_{m=m_{\alpha,n}}^{M}\frac{\exp(-2\pi mn)}{(2\pi m+2\pi w+4\pi i b n)^{\alpha+1}}\,. 
\end{split}
\end{equation}

\noindent
We claim expression (\ref{eq:aux4}) can be evaluated to within $\pm \,100\,\nu(K,j,\epsilon)\,\epsilon$ using $\le 1000\,\nu(K,j,\epsilon)^2$ operations on numbers of \mbox{$100\,\nu(K,j,\epsilon)^2$} bits. To see why, notice if $n\ne 0$, the series over $m$ can be truncated after $L:=L(K,j,\epsilon)$ terms, with a truncation error $\le 10\,(\alpha+1)\,\exp(-2\pi n (\alpha+L)) \le 10\,\epsilon$, where we used the facts $\alpha^v\le (m_{\alpha,n}+w+2bn)^v$, which holds by construction, and $\alpha !/(\alpha +1-v)!\le \alpha^v$. Once truncated, the series (\ref{eq:aux4}) can be evaluated directly in $\le 100L(K,j,\epsilon)$ operations.  If $n=0$, the series (\ref{eq:aux4}) is equal to

\begin{equation} \label{eq:aux401}
-\sum_{v=1}^{\alpha+1} \frac{\alpha !}{(\alpha+1-v)!}\sum_{m=m_{\alpha,n}}^M \frac{\exp(-2\pi m-2\pi w))}{(2\pi m+2\pi w)^v} +\, \alpha !\sum_{m=m_{\alpha,n}}^{M}\frac{1}{(2\pi m+2\pi w)^{\alpha+1}}\,. 
\end{equation}

\noindent
Since the terms in the first series over $m$ in (\ref{eq:aux401}) decline exponentially fast with $m$ (due the the decay provided by the term $e^{-2\pi m}$), it can be truncated early, after $L:=L(K,j,\epsilon)$ terms, with truncation error $\le 10\, \epsilon$. The truncated series can then be evaluated directly. As for the second series in (\ref{eq:aux401}), it can be calculated efficiently using the Euler-Maclaurin summation formula; specifically, the initial sum $\sum_{m_{\alpha,n}\le m <10 (m_{\alpha,n}+L)}$, which consists of \mbox{$\le 10(m_{\alpha,n}+L)\le 100\,\nu(K,j,\epsilon)$} terms, is evaluated directly, while the tail sum $\sum_{10(m_{\alpha,n}+L)\le m \le M}$ is evaluated to within $\pm \,10\,\nu(K,j,\epsilon)\,\epsilon$ using an Euler-Maclaurin formula like (\ref{eq:ems}) at a cost of $\le 100\,\nu(K,j,\epsilon)^2$ operations on numbers of \mbox{$\le 100\, \nu(K,j,\epsilon)^2$} bits say. 

It remains to deal with the subsum $\sum_{1\le m<m_{\alpha,n}}$ from (\ref{eq:add2}). Since this subsum consists of $<m_{\alpha,n}=2L+j\le 10\,\nu(K,j,\epsilon)$ terms, it suffices to show how to deal with a single term there, which is essentially an integral of the form 

\begin{equation} \label{eq:add201}
\int_0^{1} t^\alpha \exp \left(-2\pi (m+w +2i b n )t \right)\,dt\,,\qquad 1\le m<m_{\alpha,n}\,.
\end{equation}

\noindent
To do so, we apply the change of variable $t\leftarrow \lceil m+w+2bn\rceil\,t$ to (\ref{eq:add201}) to reduce it to a sum of the $\lceil m+w+2bn\rceil\le 10\,\nu(K,j,\epsilon)$ integrals

\begin{equation}\label{eq:anintg3}
\frac{1}{\lceil m+w+2bn\rceil^{\alpha+1}}\int_l^{l+1} t^{\alpha} \exp \left(-2\pi (m+w+2i b n)t/\lceil m+w+2bn\rceil \right)\,dt\,,
\end{equation}

\noindent
where $0 \le l\le \lceil m+w+2bn\rceil-1$ is an integer. The integrals (\ref{eq:anintg3}) are straightforward to evaluate: one makes the change of variable $t\leftarrow t-l$, then uses Taylor expansions to break down the integrand into a polynomial in $t$ of degree $2L+\alpha$ say, plus an error of size $O(\epsilon)$, and finally one integrates each term explicitly (note each term is just a monomial $z_d t^d$ for some integer $0\le d\le 2L+\alpha$, and some quickly computable coefficient $z_d$). 
\end{proof}

\begin{lem} \label{lem:a2}
There are absolute constants $\kappa_5$, $\kappa_6$, $A_7$, $A_8$, and $A_9$, such that for any positive $\epsilon <e^{-1}$, any integer $0\le j$, any integer $10\,\nu(K,j,\epsilon)^2<K$ say, any $0\le b\le 1$ satisfying $1\le 2bK$ say, and any $0\le w\le 1$ say, each of the integrals

\begin{displaymath}
\begin{array}{lll}
&\tilde{I}_{\overline{C_1}}(K,j;w,b), & \tilde{I}_{C_7}(K,j;w,b)\,, \nonumber \\
&\tilde{I}_{C_9}(K,j;w,b), &  \tilde{I}_{C_9}(K,j;w-iw+2bK +i2bK,-2ib)\,,  \nonumber 
\end{array}
\end{displaymath}

\vspace{10pt}

\noindent
can be evaluated to within $\pm \,A_7\,\nu(K,j,\epsilon)^{\kappa_{5}} \epsilon$ using $\le A_8\,\nu(K,j,\epsilon)^{\kappa_{6}}$ operations on numbers of $\le A_9\,\nu(K,l,\epsilon)^2$ bits. Moreover, under the same assumptions on $K$, $j$, and $b$, as above, except $b$ need not satisfy the condition  $1\le 2bK$, and for any $-1\le a\le 1$ say, the integral $I_{C_0}(K,j;a,b)$ can be evaluated with the same accuracy and efficiency as the above four integrals. 
\end{lem}

\begin{proof}
We show how to compute $\tilde{I}_{\overline{C_1}}(K,j;w,b)$  first. We have

\begin{equation}
\tilde{I}_{\overline{C_1}}(K,j;w,b)=  c_8 e^{-2\pi w K}\sum_{l=0}^j \binom{j}{l}\frac{(-i)^l}{K^l} \int_0^K t^l \exp \left( 2\pi iw t -4\pi bK t +2\pi i bt^2 \right) \, dt\,, \nonumber
\end{equation}

\noindent
where $c_8:=c_8(b,K)=-ie^{-2\pi i b K^2}$. Since $2bK\ge 1$ by hypothesis, we can truncate the interval of integration above at $L:=L(K,j,\epsilon)=\lceil \nu(K,j,\epsilon) \rceil$, which reduces our task to evaluating $(j+1) L$ integrals of the form

\begin{equation} \label{eq:la21}
\frac{1}{L^l}\int_n^{n+1} t^l \exp \left( 2\pi iw t -4\pi bK t +2\pi i bt^2 \right) \, dt\,,
\end{equation}

\noindent
for integers $0\le l \le j$ and $0\le n\le L-1$. To evaluate (\ref{eq:la21}), substitute $t\leftarrow t-n$, then eliminate the quadratic term $\exp(2\pi i bt^2)$ using Taylor expansion. This results in a linear combination, with quickly computable coefficients each of size $O(1)$, of, say, $3L$ integrals of the form

\begin{equation}\label{eq:anintg}
 \int_0^1 t^{\alpha} \exp \left( 2\pi i\eta t  \right) \, dt\,, 
\end{equation}

\noindent
where $\eta:=\eta_{n,w,b,K}=w+2bn+2ibK$ and $0\le \alpha <3L$ an integer. The integrals (\ref{eq:anintg}) are easily-calculable: if $\alpha<|w+2bn+2ibK|$, we evaluate (\ref{eq:anintg}) explicitly as was done in (\ref{eq:anintg2}), and if $|w+2bn+2ibK|\le \alpha$, we follow similar techniques to those used to arrive at expression (\ref{eq:anintg3}) earlier. The evaluation of $\tilde{I}_{C_9}(K,l;w-iw+2bK+i2bK,-2ib)$ is completely similar to $\tilde{I}_{\overline{C_1}}(K,j;w,b)$, already considered. 

We move on to $\tilde{I}_{C_7}(K,j;w,b)$. We have by definition

\begin{equation}
\tilde{I}_{C_7}(K,j;w,b) =\frac{c_9}{K^j} \int_0^{\sqrt{2} K} t^j \exp \left(-\sqrt{2}\pi w t +\sqrt{2}\pi iw t-2\pi b t^2 \right) \, dt\,, \nonumber
\end{equation}

\noindent
where $c_9:=c_9(j)=\exp\left(-(j+1)\pi i/4\right)$. The change of variable $t\leftarrow \sqrt{b} \,t $ yields

\begin{equation}
\tilde{I}_{C_7}(K,j;w,b) = \frac{c_9 }{ b^{(j+1)/2} K^j} \int_0^{\sqrt{2b} K} t^j \exp \left(-2\pi \frac{w}{\sqrt{2b}} t +2\pi i\frac{w}{\sqrt{2b}} t-2\pi t^2 \right) \, dt\,. \nonumber
\end{equation}

\noindent
So, truncating the interval of integration at $\lceil \sqrt{L} \rceil$ reduces the problem to evaluating 

\begin{equation}
\frac{c_9}{b^{(j+1)/2} K^j} \int_n^{n+1} t^j \exp \left(-2\pi \frac{w}{\sqrt{2b}} t +2\pi i\frac{w}{\sqrt{2b}} t-2\pi t^2 \right) \, dt\,, 
\end{equation}

\noindent
for integers $0\le n<\lceil \sqrt{L}\rceil$. The integrals are handled as follows: substitute $t\leftarrow t-n$, then eliminate the quadratic term using Taylor expansions, this results in integrals similar to (\ref{eq:anintg}), which we already know how to handle. 

Next, we consider $\tilde{I}_{C_9}(K,j;w,b)$. If $w=0$, this integral is quickly calculable via the self-similarity formula (\ref{eq:ssg}), or some variation of it. So we may assume $w>0$. Since  

\begin{equation} \label{eq:est66}
\left| \frac{1}{K^{j}} \int_{0}^{T} (T-it)^j \exp\left(- 2\pi w(T-it)-2\pi i b (T-it)^2 \right)\,dt \right|  \to_{T\to \infty} 0\,.
\end{equation}

\noindent
then by Cauchy's theorem, we may replace $C_9$ by $e^{-\pi i /4} C_9$ in $\tilde{I}_{C_9}(K,j;w,b)$. (We remark if $j=0$, then (\ref{eq:est66}) holds uniformly in \mbox{$0\le \omega\le 1$}. This observation is used in the proof of lemma~\ref{lem:next1} later.) Combined with a straightforward estimate, this yields

\begin{equation} \label{eq:est67}
\tilde{I}_{e^{-\pi i /4} C_9}(K,j;w,b)=\tilde{I}_{C_7}(K,j;w,b)+O(e^{-K})\,, 
\end{equation}

\noindent
which we have already shown how to compute.

Last, we consider the integral $I_{C_0}(K,j;a,b)$. This may contain a critical point or it may not according to whether $-a/(2b)\in [0,K]$ or not. We supplied methods to deal with  these  possibilities in \textsection{3.1.1} and \textsection{3.1.2} respectively, provided $1\le 2bK$. But the same methods still apply as long as $b$ is not too small, say $1<bK^2$. If not, say $b<1/K^2$, then computing $I_{C_0}(K,j;a,b)$ is straightforward anyway because one can apply Taylor expansions to the quadratic factor $\exp(2\pi i b t^2)$ in $I_{C_0}(K,j;a,b)$ to reduce it to a polynomial in $t$ of degree $2L$ say, plus an error of size $O(\epsilon)$, which, on applying the change of variable $t\leftarrow t/K$, yields an integral similar to (\ref{eq:anintg}), which we have already shown how to handle. 

\end{proof}

\begin{lem} \label{lem:a3}
For any integer $K>0$, any integer $j\ge 0$, any integer $m$, any $a\in \mathbb{R}$, and any $b>0$ such that $q:=\lfloor a+2bK \rfloor$ is not zero, we have 

\begin{equation} \label{eq:expform}
I_{C_8}(K,j;a-m,b)= \exp\left(\frac{2\pi i a}{2b}m-\frac{2\pi i}{4b}m^2\right)  \sum_{s=0}^j \frac{w_{s,j,a,b,K}\, m^s}{q^s}\,, 
\end{equation}

\begin{equation} \label{eq:wsj}
\begin{split}
w_{s,j,a,b,K}=q^s & \frac{j!\, \sqrt{2\pi} e^{\pi i /4}e^{(j-s)3\pi i/4} e^{-i \pi a^2/(2b)}}{\,2^{j/2}s! (2\sqrt{b\pi})^{j+1}K^j}\left(\sqrt{\frac{2\pi}{b}}\right)^{s}  \\
& \times\, \sum_{l=0}^{j-s}   \frac{\delta_{(j-s-l) \bmod{2}} (-1)^{(j+l-s)/2}}{l!\frac{j-s-l}{2}!}\left(ae^{-3\pi i /4} \sqrt{\frac{2\pi}{b}}\right)^{l}\,.
\end{split}
\end{equation}

\noindent
We remark (\ref{eq:expform}) is what one would expect; it is also essentially independent of $K$. The normalization by $q^s$, as well as the shifting by $m$, in the statement of the lemma is done because it is convenient in the context of our proof of Theorem~\ref{thm: mainthm} in \textsection{3} and \textsection{4}. 
\end{lem}

\begin{proof}
This follows from well-known properties of the Hermite polynomials; see~\cite{Is}.
\end{proof}

\begin{lem} \label{lem:a4}
For any $\epsilon \in (0,e^{-1})$, any $a \in [0,1]$, any $b \in [0,1]$, any integer $j\ge 0$, any positive integer $K>\Lambda(K,j,\epsilon)$, any integer $0\le s \le j$, let $w_{s,m,a,b,K}$ be defined as in (\ref{eq:wsj}), then assuming $\lceil a \rceil <\lfloor a+2bK \rfloor$, we have 

\begin{equation}\label{eq:wbbbd}
\sum_{m=s}^j |w_{s,m,a,b,K}| \le \frac{e}{ \sqrt{2b}} \left( 1+\frac{1}{2bK}\right)^j \sum_{g=0}^{j}  \left(\frac{j}{2bK}\right)^g \,.
\end{equation}

\noindent
If in addition $2bK \le  4\,\nu(K,j,\epsilon)^3$ say, then  $\sum_{m=s}^j |w_{s,m,a,b,K}| \le (2b)^{-1/2} (j+1) 4^{j+2}$.
\end{lem}

\begin{proof}
From formula (\ref{eq:wsj}), and the bounds $b\in [0,1]$ and $s\in [0, j]$, we obtain 

\begin{equation}
\begin{split}
\sum_{m=s}^j |w_{s,m,a,b,K}|& \le  \frac{(\lfloor a+2bK\rfloor)^s}{(2bK)^{s}}\frac{1}{\sqrt{2b}} \sum_{m=0}^{j-s}  \frac{(m+s)^m}{\,\left(\sqrt{2\pi}\right)^{m}(2bK)^{m}}  \sum_{\substack{0\le l\le m\\m-l \textrm{ even}}} \frac{\left(\sqrt{2\pi} a\right)^l b^{(m-l)/2}}{l! \,\frac{m-l}{2}!} \\ 
&\le \left( 1+\frac{a}{2bK}\right)^j \frac{1}{\sqrt{2b}}  \left[\sum_{g=0}^{j}  \left(\frac{j}{2bK}\right)^g\right]\,e^a\,. \nonumber
\end{split}
\end{equation}

\noindent
The bound (\ref{eq:wbbbd}) now follows because $a\in [0,1]$ by hypothesis. To prove the last part of the lemma, notice if $2bK \le 4\, \nu(K,j,\epsilon)^3$, then since $\Lambda(K,j,\epsilon)<K$, it follows $b<1/(2j+2)^2$.  Also, the assumption $\lceil a\rceil <\lfloor a+2bK\rfloor$ implies $1/(2K)\le b$. Therefore, by the definition (\ref{eq:wsj}), and a direct calculation,

\begin{equation}
\sum_{m=s}^j |w_{s,m,a,b,K}| \le  \frac{2q^s}{(2bK)^{s} \sqrt{2b}} \sum_{m=0}^{j-s}  \frac{(m+s)!}{s!\,m!\, (2bK)^{m}}\le  \frac{(j+1) 4^{j+2}}{\sqrt{2b}}\,.  
\end{equation}

\end{proof}

\begin{lem} \label{lem:a5}
For any integer $j\ge 0$, any integer $m\ge 0$, any integer $K>0$, and any real numbers $a$ and $b$, the function $f_{K,j,a,b}(x):=\frac{x^{j}}{K^{j}}\exp(2\pi i a x+2\pi i b x^2)$ satisfies

\begin{equation}
\max_{0\le x \le K}|f_{K,j,a,b}^{(m)}(x)| \le (2\pi (|a|+|2bK|)+(m+j)/K)^{m}\,. 
\end{equation}
\end{lem}

\begin{proof}
$f^{(m)}_{K,j,a,b}(x)= P_{m,K,j,a,b}(x)\exp(2\pi i a x+2\pi i b x^2)$ where $P_{m,K,j,a,b}(x)$ is a polynomial in $x$ of degree $m+j$. So $P_{m,K,j,a,b}(x):=\sum_{l=0}^{m+j} d_{l,m,K,j,a,b} \, x^l$ for some coefficients $d_{l,m,K,j,a,b}$ defined by the recursion

\begin{equation}\label{eq:recpm}
P_{m+1,K,j,a,b}(x)=2\pi i (a+2bx) P_{m,K,j,a,b}(x)+P'_{m,K,j,a,b}(x)\,,
\end{equation}

\noindent
where $P_{0,K,j,a,b}(x):=x^j/K^j$ and $P'_{m,K,j,a,b}(x)$ denotes the derivative of $P_{m,K,j,a,b}(x)$ with respect to $x$. Notice $|f_{K,j,a,b}^{(m)}(x)|=|P_{m,K,j,a,b}(x)|$. Define $|P_{m,K,j,a,b}(x)|_1:=\sum_{l=0}^{m+j} |d_{l,m,K,j,a,b} x^l|$, and notice $|P(x)|\le |P(x)|_1$. By induction on $m$, suppose 

\begin{equation}\label{eq:anbbbd}
\max_{0\le x\le K}|P_{m,K,j,a,b}(x)|_1 \le (2\pi (|a|+|2bK|)+(m+j)/K)^{m}\,.
\end{equation}

\noindent
Clearly, (\ref{eq:anbbbd}) holds when $m=0$, and it is straightforward to verify 

\begin{equation}\label{eq:recpm2}
\max_{0\le x \le K} |P'_{m,K,j,a,b}(x)|_1 \le \frac{m+j}{K}\, \max_{0\le x\le K} |P_{m,K,j,a,b}(x)|_1\,. 
\end{equation}

\noindent
On combining relations (\ref{eq:recpm}) and (\ref{eq:recpm2}), we obtain

\begin{equation}
\begin{split}
\max_{0\le x\le K} |P_{m+1,K,j,a,b}(x)|_1 & \le \max_{0\le x \le K} |2\pi i (a+2bx) P_{m,K,j,a,b}(x)|_1  \\
& +\max_{0\le x \le K}|P'_{m,K,j,a,b}(x) |_1  \\
& \le (2\pi (|a|+|2bK|)+(m+1+j)/K)^{m+1}\,,
\end{split}
\end{equation}

\noindent
as required. Notice the inductive proof naturally gives a method to compute the polynomials $P_{m,K,j,a,b}(x)$.
\end{proof}

\begin{lem} \label{lem:next1}
Let $\epsilon \in (0,e^{-1})$, $a \in [0,2]$, $b\in [0,1/4]$, and $K>0$ an integer. Define $\nu(K,\epsilon):=\log (K/\epsilon)$, $M:=M(K,\epsilon)=\lceil K^3 e^{\nu(K,\epsilon)}\rceil$, $F(K;a,b):=F(K,0;a,b)$,  $p_a=\lceil a\rceil$, $q_a:=q_{a,b,K}= \lfloor a+2bK\rfloor$, $p_{1,a}:=p_{1,a,b,K}= q_{a,b,K}-p_{a,b,K}$, $\omega_a:=\omega_{a,b,K}=\{a+2bK\}$, and $\omega_{1,a}=p_a -a$. Let $\delta_n$ denote the function which is 1 for $n=0$, and 0 otherwise, and let $J(.)$ and $\tilde{I}_{C_.}(.)$ be as defined in \textsection{2}. Then for any tuple $(\alpha,a,b) \in [-1,1]\times [0,2]\times [0,1/4]$ such that $p_{a+\alpha x} < q_{a+\alpha x}$ and $a+\alpha x \in (0,2)$ for all $x\in [-1/4,1/4]$, we have

\begin{equation}
\begin{split}
F(K;a+\alpha x,b)= & e^{\pi i /4 -\pi i (a+\alpha x)^2/(2b)} \, F\left(\lfloor 2bK\rfloor;\frac{a+\alpha x}{2b},-\frac{1}{4b}\right)\\
& +R_{M}(K,a+\alpha x,b)  +O(K^{-2}\epsilon +e^{-K})\,,
\end{split} 
\end{equation}

\vspace{10pt}

\noindent
where $x$ is any number in $[-1/4,1/4]$, and $R_{M}(K,a+\alpha x,b)$ is a linear combination of the constant function 1, and the following eighteen functions: 

\begin{displaymath}
\begin{array}{ll}
J(K;M,2bK-\omega_{a+\alpha x},b)\,, &  e^{2\pi i \alpha x K}\, J(K;M,2bK-\omega_{1,a+\alpha x},b)\,, \\
J(K;p_{1,a+\alpha x},\omega_{1,a+\alpha x},b)\,, & e^{2\pi i \alpha x K}\, J(K;p_{1,a+\alpha x},\omega_{a+\alpha x},b)\,,\\ 
J(K;M,1-\omega_{1,a+\alpha x},b)\,, &  e^{2\pi i \alpha x K}\, J(K;M,1-\omega_{a+\alpha x},b)\,, \\
\tilde{I}_{C_7}(K;1-\omega_{1,a+\alpha x},b)\,, &  e^{2\pi i \alpha x K} \, \tilde{I}_{C_7}(K;1-\omega_{a+\alpha x},b)\,,\\
\tilde{I}_{C_7}(K;\omega_{1,a+\alpha x},b)\,, & e^{2\pi i \alpha x K} \, \tilde{I}_{C_7}(K;\omega_{a+\alpha x},b)\,, 
\end{array}
\end{displaymath}

\vspace{10pt}

\begin{displaymath}
\begin{array}{ll}
\frac{1}{\sqrt{2b}}\,e^{-\pi i (a+\alpha x)^2/(2b)}\,, & e^{2\pi i \alpha x K -2\pi \omega_{a+\alpha x} K}\, \tilde{I}_{C_0}(K; e^{\pi i /4}(-i\omega_{a+\alpha x} + 2bK),-ib)\,,\\
e^{2\pi i \alpha x K}\,, & e^{2\pi i \alpha x K-2\pi \omega_{a+\alpha x} K} \, \tilde{I}_{C_0}(K; -i\omega_{a+\alpha x} +2bK,-b)\,.
\end{array}
\end{displaymath}

\vspace{10pt}

\begin{displaymath}
\begin{array}{ll}
c_{1,a+\alpha x}\, e^{2\pi i \alpha x/(2b) -\pi i (a+\alpha x)^2/(2b)}\,, & c_{2,a+\alpha x} \, e^{2\pi i \alpha x(K^*+1) /(2b)-\pi i (a+\alpha x)^2/(2b)}\,, \\
c_{3,a+\alpha x} \, e^{2\pi i \alpha x(K^*+1) /(2b) -\pi i (a+\alpha x)^2/(2b)}\,, & c_{3,a+\alpha x}\, e^{2\pi i \alpha x(K^*+2) /(2b)-\pi i (a+\alpha x)^2/(2b)}\,,
\end{array}
\end{displaymath}

\vspace{10pt}

\noindent
where $c_{1,a}=\delta_{2-p_a}$, $c_{2,a}:=c_{2,a,b,K}=\delta_{q_{a,b,K} - K^*_{b,K}-1}$, and $c_{3,a}:=c_{3,a,b,K} = \delta_{q_{a,b,K} - K^*_{b,K}-2}$. The coefficients in the linear combination can all be computed to within $\pm\, \epsilon/K^2$ say using $O(\nu(K,\epsilon))$ operations on numbers of $O(\nu(K,\epsilon))$ bits, are bounded by $O(1)$, and do not depend on $x$. Implicit asymptotic constants are absolute.
\end{lem}

\begin{proof}
This follows directly from formulas (\ref{eq:sum1}), (\ref{eq:sum2}), and (\ref{eq:sum3}), the method of proof of lemmas~\ref{lem:a2} and~\ref{lem:a3}, the remarks following formulas (\ref{eq:est63}) and (\ref{eq:est66}), and some routine calculations and estimates. The conditions $p_{a+\alpha x}< q_{a+\alpha x}$ and $a+\alpha x\in (0,2)$ for all $x\in [-1/4,1/4]$, which are stated in the lemma, are not essential but they help simplify the presentation of lemma~\ref{lem:next2} next.  
\end{proof}

\begin{lem} \label{lem:next2}
Let $\epsilon \in (0,e^{-1})$, $K>\Lambda(K,\epsilon):=1000\,\nu(K,\epsilon)^6$ say, $K$ an integer, and $(\alpha,a,b)\in [-1/\Lambda(K,\epsilon),1/\Lambda(K,\epsilon)]\times [0,2]\times [0,1/4]$. Let $[w,z)\subset [-1/4,1/4]$ be any subinterval such that $p_{a+\alpha x}$ and $q_{a+\alpha x}$ are constant over $x\in [w,z)$, $p_{a+\alpha x}<q_{a+\alpha x}$ for all $x\in [w,z)$, and $a+\alpha x \in (0,2)$ for all $x\in [w,z)$. Last, let $l$ and $m$ denote any integers satisfying $ m,l \in [0, 1000\,\nu(K,\epsilon)]$ say. Then for any $x\in [w,z)$,  each of the eighteen functions listed in lemma~\ref{lem:next1} can be written as a linear combination of the functions

\begin{displaymath}
 x^m\,, \qquad  x^m \exp\left(2\pi i \alpha x K\right)\,, \qquad  \exp\left(2\pi i \alpha x P /(2b)-2\pi i \alpha^2 x^2/(4b)\right)\,,  \\
\end{displaymath} 

\vspace{10pt}

\noindent
where $P\in \{-1,0,K^*,K^*+1 \}$, and the functions

\begin{displaymath}
\begin{split}
&\exp\left( 2\pi i \,\omega_{a+\alpha x}\,N - 2\pi (1-i)m\, \frac{\omega_{a+\alpha x}}{\sqrt{2b}}\right)\times \\
& \qquad\qquad\qquad\qquad \int_0^1 t^l \exp\left( - 2\pi (1-i)\, \frac{\omega_{a+\alpha x}}{\sqrt{2b}}t-2 \pi m t\right)\,dt\,, 
\end{split}
\end{displaymath}

\vspace{10pt}

\noindent
where $N\in \{0,K\}$, and the functions 

\begin{displaymath}
(\omega_{a+\alpha x})^m\, \exp\left(2\pi i \, \omega_{a+\alpha x}\, L -  2\pi \,\omega_{a+\alpha x}\,R \right)\,,  \\
\end{displaymath}

\vspace{10pt}

\noindent
where $ L, R \in [K,  K+ 1000 \,\nu(K,\epsilon)]$ say, as well as functions of the same form, but with $\omega_{a+\alpha x}$ possibly replaced  by $1-\omega_{a+\alpha x}$ or $\omega_{1,a+\alpha x}$ or $1-\omega_{1,a+\alpha x}$, plus an error term bounded by $O(\Lambda(K,\epsilon) K^{-2} \epsilon)$. The length of the linear combination is $O(\nu(K,\epsilon))$ terms. The coefficients in the linear combinations can all be computed to within $\pm\, \epsilon/K^2$ using $O(\Lambda(K,\epsilon))$ operations on numbers of $O(\nu(K,\epsilon)^2)$ bits, are bounded by $O(K)$, and are independent of $x$.  Implicit Big-$O$ constants are absolute.
\end{lem}

\begin{proof}
This follows from lemma~\ref{lem:next1}, the proofs of lemmas~\ref{lem:a1} and~\ref{lem:a2}, the assumption that $p_{a+\alpha x}$ and $q_{a+\alpha x}$ are constant over $x\in [w,z)$, and some routine calculations.
\end{proof}

\textbf{Acknowledgment.} I would like thank my PhD thesis advisor Andrew Odlyzko. Without his help and comments this paper would not have been possible. I would like to thank Jonathan Bober, Dennis Hejhal, and Michael Rubinstein, for helpful remarks.

\end{document}